\documentclass[draft,utf8,12pt]{aclart}
\usepackage{lmodern}
\DeclareUnicodeCharacter{00A0}{\zztop\nobreakspace}

\title{Mesures de Mahler et \'equidistribution logarithmique}
\author{Antoine Chambert-Loir}
\address{Universit\'e de Rennes~1 \\ 
IRMAR (UMR 6625 du CNRS) \\
Campus de Beaulieu \\ 35042 Rennes Cedex}
\email{antoine.chambert-loir@univ-rennes1.fr}
\author{Amaury Thuillier}
\address{Universit\'e de Lyon \\
Universit\'e Lyon1 \\ CNRS, UMR 5208 Institut Camille Jordan \\
B\^atiment du Doyen Jean Braconnier \\ 43 bd. du 11 novembre 1918 \\
F - 69622 Villeurbanne Cedex} \email{thuillier@math.univ-lyon1.fr}

\def\C{{\mathbf C}}
\def\P{{\mathbf P}}
\def\Q{{\mathbf Q}}
\def\N{{\mathbf N}}
\def\R{{\mathbf R}}
\def\Z{{\mathbf Z}}
\def\sp{\operatorname{sp}}
\def\ord{\operatorname{ord}}
\def\dim{\operatorname{dim}}
\let\ra\rightarrow
\let\hra\hookrightarrow
\let\emptyset\varnothing
\let\eps\varepsilon
\let\phi\varphi
\let\bar\overline

\def\cf{\emph{cf.}\xspace}

\def\abs#1{\left|{#1}\right|}
\def\norm#1{\left\| {#1} \right\|}

\def\hc{{\widehat c}}

\def\hdiv{\mathop{\widehat{\operatorname{div}}}}
\def\hvol{\mathop{\widehat{\operatorname{vol}}}}
\def\vol{\operatorname{vol}}

\def\Gal{\operatorname{Gal}}

\def\Spec{\operatorname{Spec}}
\def\Proj{\operatorname{Proj}}
\def\div{\operatorname{div}}

\def\ddc{\mathop{\mathrm d\mathrm d^c}}
\allowdisplaybreaks[3]

\begin{document}

\begin{abstract}
Soit $X$ un sch\'ema projectif int\`egre sur un corps de nombres~$F$; 
soit $L$ un fibr\'e en droites ample sur~$X$ muni d'une m\'etrique ad\'elique
semi-positive au sens de Zhang.
Les r\'esultats principaux de cet article sont :

1) Une formule qui calcule
les hauteurs locales (relativement \`a~$L$) d'un diviseur de Cartier sur~$X$
comme des {\og mesures de Mahler\fg} g\'en\'eralis\'ees, c'est-\`a-dire
les int\'egrales de fonctions de Green pour~$D$ contre 
des mesures associ\'ees \`a~$L$;

2)  Un th\'eor\`eme d'\'equidistribution des points de {\og petite\fg}
hauteur valable pour des fonctions-test \`a singularit\'es logarithmiques
le long d'un diviseur~$D$, pourvu que la hauteur de~$D$
soit {\og minimale\fg}. Dans le contexte de la dynamique alg\'ebrique,
{\og petite\fg} signifie de hauteur tendant vers~$0$, 
et {\og minimale\fg} signifie de hauteur nulle.
\end{abstract}

\begin{altabstract}
Let $X$ be a projective integral scheme over a number field~$F$ ; 
let $L$ be a ample line bundle on~$X$ 
together with a semi-positive adelic metric in the sense of Zhang.
The main results of this article are

1) A formula which allows to compute the local heights
(relative to $L$) of a Cartier divisor~$D$ on~$X$ as 
generalized ``Mahler measures'',
\emph{i.e.},  integrals of Green functions for~$D$ 
against measures attached to~$L$;

2) A theorem of equidistribution of points of ``small'' height valid for functions with logarithmic singularities along a divisor~$D$, provided the height of~$D$ is ``minimal''. In the context of algebraic dynamics, ``small'' means of height converging to~$0$, and ``minimal'' means height~$0$.
\end{altabstract}

\maketitle

\section{Introduction}

Soit $F$ un corps de nombres, soit $X$ un $F$-sch\'ema projectif
int\`egre, et soit $L$  un fibr\'e en droites
ample sur~$X$, muni d'une m\'etrique ad\'elique semi-positive au
sens de~\cite{zhang95b}. Un tel fibr\'e en droites m\'etris\'e
\'etant donn\'e, tout sous-sch\'ema ferm\'e int\`egre~$Z$ de~$X$ poss\`ede
un degr\'e~$\deg_L(Z)$ et une hauteur $h_{\bar L}(Z)$, donn\'es par
les formules
\[ \deg_L(Z) = (c_1(L)^{\dim Z}|Z)\quad\text{et}\quad
h_{\bar L}(Z)=(\hc_1(\bar L)^{\dim Z+1}|Z). \] La premi\`ere formule
fait intervenir la premi\`ere classe de Chern $c_1(L)$ de~$L$ et la
th\'eorie de l'intersection classique (voir par exemple~\cite{fulton98});
la seconde rel\`eve de  la th\'eorie de l'intersection arithm\'etique~\cite{gillet-s90} de Gillet--Soul\'e,
$\hc_1(\bar L)$ \'etant la premi\`ere classe de Chern arithm\'etique
du fibr\'e en droites m\'etris\'e~$\bar L$. 
Ces hauteurs g\'en\'eralisent aux
sous-vari\'et\'es les hauteurs au sens de Weil. Pour tout point~$x$
de~$X(\bar\Q)$, notons en effet~$[x]$ le point ferm\'e de~$X$
correspondant ; la fonction sur~$X(\bar\Q)$ d\'efinie par $x\mapsto
h_{\bar L}([x])/\deg_{L}([x])$ est une fonction hauteur pour le
fibr\'e en droites~$L$, fonction que nous noterons encore $h_{\bar
L}$.
L'utilisation de la th\'eorie de l'intersection arithm\'etique
pour d\'efinir une notion de hauteur de sous-vari\'et\'es 
en g\'eom\'etrie diophantienne est due  \`a Faltings~\cite{faltings1991}
et le formalisme a \'et\'e d\'evelopp\'e dans l'article~\cite{bost-g-s94}
de Bost, Gillet et Soul\'e. Nous utiliserons en fait une extension
de cette th\'eorie introduite par Zhang~\cite{zhang95b}
o\`u les mod\`eles entiers sont remplac\'es par une notion convenable
de m\'etriques \`a toutes les places, finies ou infinies.

Pour toute place~$v$ de~$F$, nous noterons $\C_v$ le compl\'et\'e
d'une cl\^oture alg\'ebrique du compl\'et\'e~$F_v$ de~$F$ en la
place~$v$ et $X_v$ l'espace analytique de~$X$ en~$v$: si $v$ est
archim\'edienne, il s'agit de l'espace analytique complexe
$X(\C_v)$; si $v$ est ultram\'etrique, il s'agit de l'espace
analytique $v$-adique au sens de Berkovich~\cite{berkovich1990}
attach\'e au sch\'ema
projectif~$X_{\C_v}$. Lorsque $v$ est archim\'edienne, diverses
formules de g\'eom\'etrie d'Arakelov, dont la d\'efinition
par r\'ecurrence de la hauteur d'une sous-vari\'et\'e
(cf. la formule en haut de la page~554 de~\cite{faltings1991}
ou la formule~(3.2.2) de~\cite{bost-g-s94})
font intervenir une mesure $\mu_{\bar L,v}$ sur~$X_v$
d\'eduite de la m\'etrique $v$-adique de~$L$ :
si la m\'etrique de~$\bar L$ est lisse, cette mesure est donn\'ee
par $c_1(\bar L)^{\dim X}/\deg_L(X)$, o\`u $c_1(\bar L)$ d\'esigne la forme
de courbure de~$\bar L$.
Dans le cas o\`u $v$ est ultram\'etrique, la construction de mesures
analogues, not\'ees $c_1(\bar L)_v^{\dim X}$ et~$\mu_{\bar L,v}$,
est faite dans l'article~\cite{chambert-loir2006} du premier auteur.
(Modulo la r\'esolution des singularit\'es
\`a laquelle elle est conditionn\'ee,
la th\'eorie d'Arakelov non archim\'edienne de~\cite{bloch-g-s95b}
permet en principe d'\'etudier le {\og courant\fg} de courbure
$c_1(\bar L)$ et ses puissances; le lien pr\'ecis
avec la th\'eorie de~\cite{chambert-loir2006}
m\'eriterait d'\^etre explicit\'e.)

Ces mesures donnent lieu \`a un th\'eor\`eme d'\'equidistribution
des points de petite hauteur.
Si $\{x_1,\dots,x_d\}$ est l'orbite sous $\Gal(\bar F/F)$ d'un point
$x \in X(\bar F)$,
on notera $\mu_{x,v}$ la mesure de probabilit\'e sur~$X_v$, somme des masses
de Dirac en les~$x_i$ divis\'ee par~$d$.

\begin{theo}\label{theo.equi}
Soit $(x_n)$ une suite \emph{g\'en\'erique} de points de~$X(\bar F)$, 
au sens o\`u chaque sous-vari\'et\'e 
$Z\subsetneq X \otimes_{F} \bar F$ 
ne contient qu'un nombre fini de termes de la suite.
\begin{enumerate}\def\theenumi{\alph{enumi}}\def\labelenumi{\theenumi)}
\item On a $\displaystyle \liminf_n h_{\bar L}(x_n) \geq \frac{h_{\bar L}(X)}{(\dim X+1)\deg_L(X)}$;
\item si de plus la suite $(h_{\bar L}(x_n))$ converge vers $h_{\bar L}(X)/(\dim X+1)\deg_L(X)$,
alors, pour toute place~$v$ de~$F$,
la suite $(\mu_{x_n,v})$ de mesures sur~$X_v$
converge vers la mesure $\mu_{L,v}$.
\end{enumerate}
\end{theo}
L'hypoth\`ese que la suite~$(x_n)$ est g\'en\'erique signifie
que la suite~$([x_n])$ des points ferm\'es correspondants du sch\'ema~$X$
converge vers le point g\'en\'erique de~$X$.
Elle est \'evidemment essentielle comme le montre l'exemple
d'une suite dont une infinit\'e de termes sont \'egaux \`a un point alg\'ebrique~$x$
tel que $h_{\bar L}(x)$ ne v\'erifie pas cette in\'egalit\'e.

L'in\'egalit\'e du~\emph a) est due \`a Zhang~\cite{zhang95b}; c'est
une cons\'equence du th\'eor\`eme de Hilbert-Samuel
arithm\'etique~\cite{gillet-s88,abbes-b95,zhang95b}. Lorsque $X$ est
une courbe, l'assertion \emph b) est un th\'eor\`eme
d'Autissier~\cite{autissier2001b} lorsque $v$ est archim\'edienne et
du premier auteur~\cite{chambert-loir2006} lorsque $v$ est
ultram\'etrique. (Voir aussi~\cite{rumely99,baker-r2006,favre-rl2006}
pour des cas particuliers lorsque $X$ est la droite projective,
notamment dans le contexte des capacit\'es ou des syst\`emes dynamiques.)
La g\'en\'eralisation au cas des vari\'et\'es de
dimension arbitraire est due \`a Yuan~\cite{yuan2006}, comme
corollaire de sa th\'eorie des fibr\'es en droites
arithm\'etiquement gros. La principe de la d\'emonstration
a \'et\'e introduit par Szpiro, Ullmo et Zhang
dans~\cite{szpiro-u-z97}
et consiste \`a appliquer l'in\'egalit\'e du~\emph a) 
\`a des petites variations du fibr\'e en droites m\'etris\'e~$\overline L$.

Sous les hypoth\`eses du th\'eor\`eme,
la convergence des mesures signifie que
pour toute place~$v$
de~$F$ et toute fonction continue $\phi$ sur~$X_v$,
la suite $\int_{X_v} \phi \, \mathrm d\mu_{x_n,v}$ converge
vers $\int_{X_v}\phi\, \mathrm d\mu_{\bar L,v}$.
Dans cet article, nous voulons montrer un exemple
o\`u cette convergence a lieu pour des fonctions $\phi$
\`a p\^oles logarithmiques le long d'un diviseur.

Si $D$ est un diviseur (de Cartier) effectif sur~$X$, $v$ une place de~$F$,
nous dirons qu'une fonction $\phi\colon X_v\ra\R\cup\{\pm\infty\}$
est \emph{\`a singularit\'es au plus logarithmiques le long de~$D$}
si elle est continue hors de~$D_v$
et si, pour tout point $x$ de~$X_v$, il existe  un voisinage~$U$ de~$x$
dans~$X_v$, une \'equation $f_U$ de~$D_v$ sur~$U$ et un nombre
r\'eel $c_U$ tel que l'on ait $\abs{\phi}\leq c_U \log\abs{f_U}^{-1}$
sur~$U$.

\begin{theo}\label{theo.equilog}
Soit $X$ un $F$-sch\'ema projectif int\`egre;
notons~$d$ sa dimension.
Soit~$\bar L$ un fibr\'e en droites ample muni d'une
m\'etrique ad\'elique semi-positive et consid\'erons une suite
g\'en\'erique~$(x_n)$ de points de~$X(\bar F)$ telle que $h_{\bar
L}(x_n)$ converge vers \mbox{$\displaystyle\frac{h_{\bar
L}(X)}{(d+1)\deg_L(X)}$}.

\begin{enumerate}\def\theenumi{\alph{enumi}}\def\labelenumi{\theenumi)}
\item Pour tout diviseur effectif~$D$ non nul sur~$X$,
on a l'in\'egalit\'e
\[  \frac{h_{\bar L}(D)}{d \deg_L (D)}
     \geq \frac{h_{\bar L}(X)}{(d+1)\deg_L(X)}. \]
\item Soit $D$ un diviseur effectif non nul pour lequel cette in\'egalit\'e soit
en fait une \'egalit\'e. Pour toute place $v$ de~$F$ et toute
fonction $\phi$ sur~$X_v$ \`a singularit\'es au plus logarithmiques
le long de~$D$, on a
\[ \int_{X_v} \phi\, \mathrm d\mu_{x_n,v} \ra \int_{X_v}\phi\,\mathrm d\mu_{\bar L,v}. \]
\end{enumerate}
\end{theo}
Observons que le membre de gauche est d\'efini d\`es que $x_n\not\in D$,
donc pour tout entier~$n$ assez grand. L'int\'egrale au membre de droite
existe aussi (\cf infra, corollaire~\ref{coro.ineg.mahler}).

\bigskip

Donnons quelques exemples, g\'en\'eralisations et commentaires.

1) Supposons que $X=\P^1$ et soit $f$ un endomorphisme de
degr\'e~$d\geq 2$ (on est en pr\'esence d'un {\og syst\`eme
dynamique alg\'ebrique\fg}). On peut alors munir~$L=\mathscr O(1)$
d'une m\'etrique ad\'elique semi-positive, unique \`a un scalaire
pr\`es, tel que $f^*\bar L$ soit isom\'etrique \`a $\bar L^{\otimes
d}$. La hauteur correspondante est appel\'ee \emph{hauteur normalis\'ee}
(par~$f$). On a alors $h_{\bar L}(X)=0$; les points de hauteur nulle
sont les points pr\'ep\'eriodiques (\emph{i.e.} ceux dont l'orbite
est finie). Dans ce cas, notre th\'eor\`eme affirme 
\emph a) la hauteur
normalis\'ee de tout diviseur effectif est positive ou nulle (un
fait classique); \emph b)
\emph{l'\'equidistribution d'une suite de
points dont la hauteur tend vers~$0$ pour une fonction \`a
singularit\'es logarithmiques en des points pr\'ep\'eriodiques.}

\medskip

2) Le cas des courbes elliptiques est analogue \`a celui des syst\`emes
dynamiques, la hauteur dynamique \'etant remplac\'ee par la hauteur
de N\'eron--Tate. Cette hauteur est positive ou nulle
et les points de hauteur nulle sont exactement les points de torsion.
 Si $X$ est une courbe elliptique,
notre th\'eor\`eme affirme ainsi
\emph{l'\'equidistribution d'une suite de points distincts dont la hauteur
de N\'eron--Tate tend vers~$0$
pour une fonction \`a singularit\'es au plus logarithmiques en
des points de torsion.}

\medskip

3) Dans le cas d'un endomorphisme $f\colon X\ra X$ d'une vari\'et\'e de
dimension sup\'erieure  munie d'un fibr\'e en droites~$L$ tel que
$f^*L\simeq L^{\otimes d}$, o\`u $d \geq 2$, un r\'esultat analogue
vaut encore: la hauteur normalis\'ee de tout diviseur est positive
ou nulle, et l'on a \'equidistribution pour des fonctions \`a
singularit\'es au plus logarithmiques le long d'un diviseur de
hauteur normalis\'ee nulle. Cependant, au moins conjecturalement,
ces diviseurs sont pr\'ep\'eriodiques et il n'en existe parfois pas.
C'est par exemple le cas lorsque $X$ est une vari\'et\'e ab\'elienne
simple et que $f$ est l'endomorphisme de multiplication par~$2$,
compte-tenu de la preuve de la conjecture de Bogomolov~\cite{ullmo98,zhang98}.
Le cas des vari\'et\'es toriques ($f$ est alors l'endomorphisme
d\'eduit de l'\'el\'evation au carr\'e sur le tore)
fournit plus d'exemples puisque les diviseurs~$D$
auquel le th\'eor\`eme s'applique sont les (adh\'erences des) sous-tores 
de codimension~$1$
et leurs translat\'es par des points de torsion.

\medskip

4) Dans le contexte dynamique des endomorphismes de~$\P^1$,
signalons aussi un th\'eor\`eme diff\'erent, mais dans le m\^eme
esprit, de Szpiro et Tucker~\cite{szpiro-t2005}
(voir aussi~\cite{everest-ward1999}). Ces auteurs
prouvent l'\'equidistribution d'une suite~$(x_n)$ de points
de~$\P^1(\bar\Q)$ tels que $h_{\bar L}(x_n)\ra 0$ vis-\`a-vis des
fonctions $\log \abs\phi$, o\`u $\phi$ est une fraction rationnelle
non nulle, pourvu que l'une des hypoth\`eses suppl\'ementaires
suivantes soit satisfaite:
\begin{itemize}
\item les points~$x_n$ sont pr\'ep\'eriodiques;
\item la suite $(f^n(x_n))$ est constante.
\end{itemize}
Leur d\'emonstration est de nature diff\'erente et repose sur
le th\'eor\`eme de Roth. Il serait int\'eressant de g\'en\'eraliser leur
th\'eor\`eme en dimension sup\'erieure.

\medskip

5) Pour toute sous-vari\'et\'e~$Y$ de~$X$,
posons $h'_{\bar L}(Y)=h_{\bar L}(Y)/(1+\dim Y)\deg_L(Y)$.
Notons que l'in\'egalit\'e~\emph a) du th\'eor\`eme~\ref{theo.equilog}
s'\'ecrit $h'_{\bar L}(D)\geq h'_{\bar L}(X)$.

\`A la suite de~\cite{autissier2006b}, il est possible
d'\'etendre le th\'eor\`eme~\ref{theo.equilog} aux suites 
de sous-vari\'et\'es de~$X$ de dimension~$p$. 
Nous d\'emontrons ainsi le th\'eor\`eme suivant.
\begin{theo}\label{theo.equilog.ssvar}
Supposons que pour toute sous-vari\'et\'e~$Z$ de dimension~$p-1$,
on ait $h'_{\bar L}(Z)\geq h'_{\bar L}(X)$.
Soit aussi $D$ un diviseur de Cartier sur~$X$ 
tel que $h'_{\bar L}(D)=h'_{\bar L}(X)$.
Alors, pour toute suite g\'en\'erique de sous-vari\'et\'es $(Y_n)$
de m\^eme dimension~$p$, telle que $h'_{\bar L}(Y_n)$ 
converge vers $h'_{\bar L}(X)$,
et toute fonction $\phi$ sur~$X_v$ \`a singularit\'es au plus logarithmiques
le long de~$D$, on a 
\[ \lim_{n\ra\infty} \frac1{\deg_L(Y_n)}
\int_{X_v}\phi c_1(\bar L)^p\delta_{Y_n,v}
=  \int_{X_v} \phi d\mu_{\bar L,v}. \]
\end{theo}
Notons que l'hypoth\`ese $h'_{\bar L}(Z)\geq h'_{\bar L}(X)$ est satisfaite
si $p=0$ (on retrouve alors le th\'eor\`eme~\ref{theo.equi})
ou si $\bar L$ est muni d'une m\'etrique ad\'elique semi-positive
associ\'ee \`a un syst\`eme dynamique.
Cela dernier cas inclut notamment les vari\'et\'es ab\'eliennes
munies de la hauteur de N\'eron--Tate (mais, comme on l'a dit
plus haut, il n'y a alors pas d'exemples int\'eressants)
et les vari\'et\'es toriques munies de leur hauteur canonique.

\medskip

6) L'in\'egalit\'e de l'assertion \emph a) est \`a mettre en regard
d'une in\'egalit\'e due \`a Zhang~\cite{zhang95}: posons
$m(\bar L)=\inf_{x\in X(\bar\Q)} h_{\bar L}(x)$
et
$e(\bar L)= \sup_{U\subset X}\inf_{x\in U(\bar \Q)} h_{\bar L}(x)$
(le {\og minimum essentiel\fg} de $h_{\bar L}$).
Alors, on a
\[ e(\bar L) \geq \frac{h_{\bar L}(X)}{(d+1)\deg_L(X)}
      \geq \frac1{d+1}(e(\bar L)+d m(\bar L)). \]
L'\'egalit\'e des deux premi\`eres quantit\'es \'equivaut \`a l'existence
d'une suite~$(x_n)$ comme dans le th\'eor\`eme.
Lorsque $X$ est une courbe,
l'assertion \emph a) affirme que
ces trois quantit\'es sont \'egales si les deux premi\`eres le sont.

\medskip

7) L'assertion d'\'equidistribution ne vaut pas
forc\'ement lorsque l'in\'egalit\'e de l'assertion \emph a) est stricte,
ainsi que le montre un contre-exemple \'el\'ementaire
d'Autissier~\cite{autissier2006}, dans lequel $X=\P^1$
et $\bar L$ est le fibr\'e en droites~$\mathscr O(1)$ muni
de la m\'etrique ad\'elique qui donne lieu \`a la hauteur de Weil.

Suivant une suggestion d'Autissier, nous montrons aussi au dernier
paragraphe de ce texte comment les m\'ethodes de cet article
permettent de produire d'autres exemples o\`u l'\'equidistribution
n'a pas lieu.

\bigskip

La d\'emonstration du th\'eor\`eme~\ref{theo.equilog}
reprend la technique qui conduit au th\'eor\`eme~\ref{theo.equi}.
Elle requiert aussi l'extension au cas de m\'etriques
\emph{admissibles}
de l'identit\'e classique en th\'eorie de l'intersection arithm\'etique
qui compare la hauteur d'une sous-vari\'et\'e
relativement \`a un fibr\'e en droites hermitien \`a celle
d'une section hyperplane
voir par exemple~\cite{bost-g-s94}, (3.2.2).
 Nous d\'emontrons ainsi  le th\'eor\`eme suivant:

\begin{theo} \label{theo.mahler}
Soit $\bar L$ et~$\bar M$ des fibr\'es en droites admissibles
sur~$X$ ;  soit aussi $s$ une section globale non nulle de~$\bar M$,
de diviseur~$D$.
Pour toute place~$v$ de~$F$,
la fonction $\log\norm{s}^{-1}_v$ est int\'egrable
pour la mesure $c_1(\bar L)^{\dim X}_v$. De plus,
\[
(\hc_1(\bar L)^{\dim X}\hc_1(\bar M)|X) = h_{\bar L}(D)
+  \sum_v \int_{X_v} \log\norm{s}^{-1}_v\,c_1(\bar L)^{\dim X}_v. \]
\end{theo}

Par d\'efinition, un fibr\'e en droites m\'etris\'e est
\emph{admissible} s'il est le quotient de deux fibr\'es en droites
munis de m\'etriques ad\'eliques semi-positives
(Zhang dit \emph{integrable} dans~\cite{zhang95b}). 

Lorsque $X= \P^n$,
ce th\'eor\`eme permet de retrouver la formule bien connue donn\'ee par
Mahler pour la hauteur d'un point de~$\P^1$, ainsi que la formule utilis\'ee
par Lawton~\cite{lawton77} dans son \'etude des hauteurs des hypersurfaces.
Supposons en effet que l'on ait $X=\P^n$ et que $D$ soit le diviseur des
z\'eros d'un polyn\^ome homog\`ene de degr\'e~$m$,
$F\in\Z[x_0,\dots,x_n]$, dont les coefficients sont premiers entre
eux. Soit $\bar L=\mathscr O(1)$, muni de la m\'etrique de Weil,
posons $\bar M=\bar L^{\otimes m}$ et soit $s$ la section de~$M$
d\'efinie par~$F$. Pour toute place~$v$ de~$\Q$ et tout point
$(x_0:\cdots:x_n)\in\P^n(\C_v)$, on a
\[ \norm{s}_v(x_0:\cdots:x_n) = \frac{\abs{F(x_0,\cdots,x_n)}_v}{\max(\abs{x_0}_v,\cdots,\abs{x_n}_v)^m}. \]
Si $v$ est une place ultram\'etrique, la mesure $c_1(\bar L)^n_v$
est la mesure de Dirac au {\og point de Gauss\fg} de $X_v$ et
l'hypoth\`ese que les coefficients de~$F$ soient des entiers
premiers entre eux signifie  que $\norm{s}_v$ vaut~$1$ en ce point.
Lorsque $v$ est la place archim\'edienne de~$\Q$, c'est la mesure de
probabilit\'e d'int\'egration  sur le polycercle:
\[ \int_{\P^n(\C)} \phi c_1(\bar L)^n_\infty
  = \frac1{(2\pi)^n} \int_{0}^{2\pi}\dots \int_0^{2\pi}
    \phi(1:e^{i\theta_1}:\cdots:e^{i\theta_n})\,d\theta_1\dots d\theta_n. \]
La formule du th\'eor\`eme~\ref{theo.mahler} devient alors
\[ m h_{\bar L}(X)= h_{\bar L}(D) + \int_{\P^n(\C)} \log\norm{s}_v^{-1} c_1(\bar L)^n_\infty. \]
En outre, $h_{\bar L}(X)=0$ (la m\'etrique de~$\bar L$ est
une m\'etrique canonique pour le syst\`eme dynamique sur~$\P^n$
donn\'e par l'\'el\'evation au carr\'e des coordonn\'ees homog\`enes), d'o\`u, finalement
\begin{align*}
 h_{\bar L}(D) & = \int_{\P^n(\C)} \log \norm{s}_v^{-1} c_1(\bar L)^n_\infty \\
&= \frac1{(2\pi)^n} \int_{0}^{2\pi}\dots \int_0^{2\pi}
    \log\abs{F(1,e^{i\theta_1},\dots,e^{i\theta_n})}\,d\theta_1\dots d\theta_n. \end{align*}
Ainsi, la hauteur d'une hypersurface de~$\P^n_\Q$ est \'egale \`a la mesure 
de Mahler du polyn\^ome homog\`ene qui la d\'efinit. Plus g\'en\'eralement,
les \emph{mesures de Mahler} du titre de cet article  d\'esignent 
les int\'egrales de fonctions de Green de diviseurs contre des mesures
de probabilit\'e.

Cette formule reliant hauteurs et nombres d'intersection arithm\'etique
est banale lorsque les m\'etriques de~$\bar L$ sont \emph{lisses}, 
c'est-\`a-dire $\mathscr C^\infty$
aux places archim\'ediennes et donn\'ees par un mod\`ele aux places
ultram\'etriques : c'est en fait une des fa\c{c}ons de d\'efinir la
hauteur d'un cycle (\cf\cite{bost-g-s94}, (3.2.2)). Elle est 
presque tautologique
lorsque $\bar M$ est le fibr\'e en droite trivial; dans ce cas,
$\log\norm{s}^{-1}$ est une fonction continue et la formule traduit
la d\'efinition des mesures $c_1(\bar L)^{\dim X}_v$. 

En fait, la difficult\'e de ce th\'eor\`eme 
repose dans l'interaction entre le processus limite qui
donne lieu aux mesures $c_1(\bar L)^d_v$ d'une part et les p\^oles
des fonctions~$\log\norm{s}^{-1}_v$ d'autre part. 
Lorsque les m\'etriques sont lisses aux places ultram\'etriques,
elle a d'abord \'et\'e d\'emontr\'ee par Maillot (\cite{maillot2000},
th\'eor\`eme 5.5.6, formule~6) comme cons\'equence
des propri\'et\'es de continuit\'e de l'op\'erateur de
Monge-Amp\`ere complexe, \'etablies par Bedford et Taylor
dans~\cite{bedford-t82} (voir aussi~\cite{demailly93}).
\`A des termes d'erreur explicites pr\`es,
cette formule figure aussi dans l'article~\cite{philippon91} de Philippon.
Lorsque $X$ est une courbe, elle a r\'ecemment fait
l'objet de plusieurs travaux ind\'ependants  : Szpiro, Tucker,
Pineiro, \cite{szpiro-t-p2005}, puis Favre et Rivera-Letelier,
\cite{favre-rl2006}, pour $X= \P^1$; Thuillier,
\cite{thuillier2005}, prop.~4.2.9 dans le cas g\'en\'eral.

\medskip

Nous d\'emontrerons au paragraphe~\ref{sec.mahler}
une version locale, valable place par place du th\'eor\`eme~\ref{theo.mahler}
et dont ce dernier se d\'eduit par sommation sur l'ensemble
des places du corps de nombres~$F$.
Nous aurons \`a faire usage d'un th\'eor\`eme d'approximation
de fonctions de Green par des fonctions lisses, th\'eor\`eme
qui fait l'objet du paragraphe~\ref{sec.approximation}.
Le th\'eor\`eme d'\'equidistribution pour les fonctions
\`a singularit\'es logarithmiques 
est d\'emontr\'e au paragraphe~\ref{sec.equilog}.
Nous avons regroup\'e dans les paragraphes~\ref{sec.rappels} et~\ref{sec.volumes}
rappels, reformulations ou compl\'ements de r\'esultats ant\'erieurs
dont nous avons besoin: le premier concerne la d\'efinition des mesures
sur les espaces analytiques ultram\'etriques, le second
formule en termes de volumes arithm\'etiques le th\'eor\`eme de Yuan~\cite{yuan2006}.
Nous concluons cet article par un \'enonc\'e de non-\'equidistribution
que nous a sugg\'er\'e P.~Autissier.

Pour cela, et pour sa lecture attentive, nous voudrions le remercier ici.
Nous sommes \'egalement reconnaissant envers S\'ebastien Boucksom
de ses remarques sur un \'etat ant\'erieur de ce travail ; nous
le remercions en particulier d'avoir attir\'e notre attention sur
une formulation erron\'ee du lemme~\ref{lemm.fatou}.
Enfin, nous sommes redevables
\`a Vincent Maillot d'avoir point\'e quelques malheurs d'expression
dans la premi\`ere version de cet article.

\section{Quelques rappels sur les fibr\'es en droites
m\'etris\'es et les mesures qui leur sont attach\'ees}
\label{sec.rappels}

Fixons une place~$v$ de~$F$. Soit $X$ un $F_v$-sch\'ema projectif
int\`egre. On note $X_v$ l'espace analytique de~$X$
sur~$\C_v$, au sens usuel si $v$ est archim\'edienne, et au sens de
Berkovich~\cite{berkovich1990} si $v$ est ultram\'etrique.

\subsection*{Mod\`eles}
Supposons que $v$ soit une place ultram\'etrique de~$F$
et notons $R$ l'anneau des entiers du corps local $F_{v}$. Par d\'efinition, un
\emph{mod\`ele} de $X$ sur $R$ est un $R$-sch\'ema $\mathscr{X}$ projectif
et plat
dont la fibre g\'en\'erique est isomorphe \`a $X$ ;
on d\'esigne alors par $\mathscr X_s$ la fibre sp\'eciale de~$\mathscr X$
et par $\sp_{\mathscr{X}} \colon X_{v} \rightarrow \mathscr{X}_s$
l'application de r\'eduction.

\begin{lemm} \label{lemm.int.cls} 
Soit $R$ un anneau de valuation
discr\`ete complet, de corps des fractions $K$, et soit $\mathcal{X}$
un $R$-sch\'ema de pr\'esentation finie et plat, de fibre g\'en\'erique
$X = \mathscr{X} \otimes_R K$.

\begin{enumerate}
\item
L'application de r\'eduction $sp_{\mathscr{X}} : X_v \rightarrow
\mathscr{X}_s$ est surjective.  
\item 
Soit $\mathscr{X}'$ la fermeture int\'egrale
de $\mathscr{X}$ dans sa fibre g\'en\'erique $X$. Le $R$-sch\'ema
$\mathscr{X}'$ est plat, et il est de pr\'esentation finie si
$\mathscr{X}$ est un sch\'ema r\'eduit.  
\end{enumerate}

Supposons de plus que $\mathscr X$ soit int\'egralement ferm\'e dans
sa fibre g\'en\'erique. Alors:
\begin{enumerate}\setcounter{enumi}{2}
\item 
Le sch\'ema $\mathscr{X}$ est  r\'eduit au voisinage de sa fibre
sp\'eciale.
\item
Pour tout
point g\'en\'erique $\xi$ de la fibre sp\'eciale $\mathcal{X}_s$,
l'anneau local $\mathcal{O}_{\mathscr{X}, \xi}$ 
est un anneau de valuation discr\`ete qui domine~$R$ ;  la semi-norme
prolongeant la valeur absolue de~$K$  et associ\'ee \`a cette  valuation
est l'unique ant\'ec\'edent du point~$\xi$ par l'application de r\'eduction
$\sp_{\mathscr{X}} \colon X_v \rightarrow \mathscr{X}_s$.
\end{enumerate}
\end{lemm}

\begin{proof}
(1) On note $\mathfrak{m}$ l'id\'eal maximal de $R$. Soit $U$ un ouvert affine de $\mathscr{X}$ et soit
$\widehat{\mathcal{A}}$ la $R$-alg\`ebre topologique
compl\'et\'ee-s\'epar\'ee de $\mathcal{A} = \mathcal{O}_{\mathscr{X}}(U)$
pour la topologie $\mathfrak{m}$-adique. 
L'anneau $B = \widehat{\mathcal{A}} \otimes_R K$ est 
une alg\`ebre strictement $K$-affino\"{\i}de.
D\'esignons par~$\tilde K$ le corps r\'esiduel de~$R$
et par~$\widetilde{B}$ 
la $\widetilde{K}$-alg\`ebre quotient de l'anneau $B^{\circ}$ des
\'el\'ements de $B$ de norme spectrale inf\'erieure \`a $1$ par
l'id\'eal $B^{\circ \circ}$ form\'es de ceux dont la norme spectrale
est strictement inf\'erieure \`a $1$.
L'application de r\'eduction $r :
\sp_{\mathscr{X}}^{-1}(U_s) = \mathcal{M}(B) \rightarrow
\Spec(\widetilde{B})$ est alors surjective 
(\cite{berkovich1990}, proposition 2.4.4).

L'anneau $B^{\circ}$ est la fermeture int\'egrale de
$\widehat{\mathcal{A}}$ dans $B$ (\cite{bosch-g-r1984}, 6.1.2 et 6.3.4). 
On dispose par cons\'equent d'un $\widetilde{K}$-morphisme \emph{surjectif}
\[ p \colon \Spec(\widetilde{B}) \rightarrow
\Spec(\widehat{\mathcal{A}} \otimes_R \widetilde{K}) = U_s \]
tel que $\sp_{U} = p \circ r$ et la surjectivit\'e de l'application
$\sp_{U} \colon \sp_{\mathscr{X}}^{-1}(U_s) \rightarrow U_s$ en d\'ecoule
imm\'ediatement.

(2) Quel que soit l'ouvert affine $U$ de $\mathscr{X}$,
la $R$-alg\`ebre plate $\mathcal{A} = \mathcal{O}_{\mathscr{X}}(U)$
s'identifie \`a un sous-anneau de~${A} = \mathcal{O}_{X}(U_K)
= \mathcal{O}_{\mathscr{X}}(U) \otimes_R K$. Par d\'efinition
de la fermeture int\'egrale, l'image r\'eciproque de~$U$ dans
le~$\mathscr X$-sch\'ema~$\mathscr X'$ est le sch\'ema affine dont
l'alg\`ebre $\mathcal{A}'$
est la fermeture int\'egrale de~$\mathcal{A}$ dans~$A$. 
Par construction, $\mathscr A'$ est sans $R$-torsion, si
bien que le $R$-sch\'ema~$\mathscr X'$ est plat.

Le $R$-sch\'ema~$\mathscr X$ est plat et excellent. Il est
alors r\'eduit si~$X$ l'est, et sa fermeture int\'egrale dans~$X$
est finie sur~$\mathscr X$, donc de pr\'esentation finie 
car $R$ est noeth\'erien.

\medskip
Supposons maintenant que le $R$-sch\'ema $\mathscr{X}$ soit int\'egralement
ferm\'e dans sa fibre g\'en\'erique.

(3)  Soit $U$ un ouvert affine
de~$\mathscr{X}$ d'alg\`ebre~$\mathscr A$; notons encore $A=\mathscr A\otimes K$.
Soit $a \in A$ un \'el\'ement nilpotent; pour
tout $t \in K$, l'\'el\'ement $t a$ de $A$ est entier sur
$\mathcal{A}$ et appartient donc \`a $\mathcal{A}$ puisque ce dernier
anneau est int\'egralement ferm\'e dans~$A$. L'\'el\'ement $a$ de $\mathcal{A}$ appartient donc \`a l'id\'eal
$\mathfrak{I} = \bigcap_{n \in \mathbb{N}} \mathfrak{m}^{n}\mathcal{A}$.
D'apr\`es le th\'eor\`eme d'intersection de Krull, il existe $b \in
\mathfrak{m}\mathcal{A}$ tel que $(1+b)\mathfrak{I} = 0$. Par
suite, l'id\'eal~$\mathfrak{I}$ et, \emph{a fortiori}, 
le nilradical de~$\mathcal{A}$ sont nuls sur l'ouvert de $\Spec(\mathcal{A})$
sur lequel $1+b$ est inversible.
Comme $b\in\mathfrak m\mathcal A$, cet ouvert contient la fibre sp\'eciale
de~$\mathscr X$.

(4) Soit $\xi$ un point g\'en\'erique de la fibre sp\'eciale
$\mathscr{X}_{s}$ de $\mathscr{X}$. L'anneau local
$\mathcal{O}_{\mathscr{X}, \xi}$ est noeth\'erien,
r\'eduit d'apr\`es l'assertion~(3); c'est un sous-anneau 
de l'anneau $\mathcal{O}_{\mathscr{X}, \xi}~\otimes_R~K$ 
qui est, lui-aussi, noeth\'erien et r\'eduit.
En vertu de la platitude de~$\mathscr X$, 
\[ \dim (\mathcal{O}_{\mathscr{X},
\xi} \otimes_R K)= \dim (\mathcal{O}_{\mathscr{X}_{s}, \xi}) \quad;\]
cette derni\`ere dimension est nulle, car $\xi$ est un point g\'en\'erique de~$\mathscr X_s$.
Par suite, l'anneau $\mathscr O_{\mathscr X,\xi}\otimes_R K$
est le produit des corps r\'esiduels de ses id\'eaux premiers minimaux,
c'est-\`a-dire le produit
des corps r\'esiduels des points g\'en\'eriques de
$\Spec(\mathcal{O}_{\mathscr{X},\xi})$ (par platitude de~$\mathscr X$ sur~$R$).
Autrement dit, $\mathscr O_{\mathscr X,\xi}\otimes K$
est l'anneau total des fractions de~$\mathscr O_{\mathscr X,\xi}$.

Par cons\'equent, le fait que $\mathscr X$ soit int\'egralement
ferm\'e dans sa fibre g\'en\'erique entra\^{\i}ne
que l'anneau local $\mathcal{O}_{\mathscr{X}, \xi}$ 
est int\'egralement ferm\'e dans son anneau total des fractions.
Comme il est de dimension~$1$, c'est un anneau de valuation
discr\`ete.

Puisque l'anneau $\mathcal{O}_{\mathscr{X}, \xi}$ domine~$R$, il existe
alors une unique valeur absolue sur~$\mathcal{O}_{\mathscr{X}, \xi}$
associ\'ee \`a la valuation $\ord_{\xi}$ qui prolonge la valeur
absolue de~$R$. Cette valeur absolue d\'efinit un point~$x$ de~$X_v$
tel que $\sp_{\mathscr{X}}(x) = \xi$. 

 
Inversement, soit $y$ un point de $X_v$ tel que $\sp_{\mathscr X}(y)=\xi$.
Le point~$y$ appartient \`a $\sp_{\mathscr{X}}^{-1}(U)$ 
et d\'efinit donc une semi-norme
multiplicative et born\'ee sur $\mathcal{O}_{\mathscr{X}}(U) \otimes_R K$ 
telle que $\abs f(y) = 1$ pour tout $f \in \mathcal{O}_{\mathscr{X}}(U)$
ne s'annulant pas au point $\xi$. Cette semi-norme se prolonge donc
\`a l'anneau local $\mathcal{O}_{\mathscr{X},\xi}$ et, d\'esignant
par $u$ une uniformisante de cet anneau de valuation discr\`ete, on a 
\[ \abs f(y) = \abs u(y)^{\ord_{\xi}(f)}\]
pour tout $f$. Par cons\'equent, $y$ est une valeur absolue sur
$\mathcal{O}_{\mathscr{X}, \xi}$ associ\'ee \`a la valuation
$\ord_{\xi}$ et prolongeant la valeur absolue de $R$, ce qui implique
$y=x$. 
\end{proof}

\subsection*{M\'etriques}
Nous aurons \`a consid\'erer divers types de fibr\'es en droites
m\'etris\'es sur~$X_v$.
Renvoyant \`a~\cite{zhang95b} ou~\cite{chambert-loir2006}
pour plus de pr\'ecisions, les d\'efinitions sont les suivantes.
Pr\'ecisons d'abord qu'il ne sera question que de m\'etriques continues.

Supposons $v$ archim\'edienne. On dit qu'un fibr\'e en droites
m\'etris\'e $\bar L$ sur $X_v$ est \emph{lisse} s'il existe, pour
tout point de $X_v$, un voisinage ouvert $U$, une immersion ferm\'ee
$U \hookrightarrow Y$ dans un espace analytique lisse et un fibr\'e
en droites $\bar{M}$ sur $Y$ muni d'une m\'etrique $C^\infty$ tels
que $\bar{L}_{|U} = i^{*}\bar{M}$. On dit que $\bar{L}$ est
\emph{semi-positif} si sa m\'etrique est limite uniforme de m\'etriques
lisses dont la forme de courbure est positive (\emph{i.e.} localement
induites par des m\'etriques $C^\infty$ dont la forme de courbure
est positive).

Supposons que $v$ soit une place ultram\'etrique. Consid\'erons un
mod\`ele~$(\mathscr X,\mathscr L)$ de~$(X,L)$, au sens o\`u
$\mathscr X$ est un mod\`ele de~$X$ sur~$R$ et o\`u $\mathscr L$ est
un fibr\'e en droites sur~$\mathscr X$ dont la restriction \`a~$X$
est identifi\'ee \`a~$L$. Ces donn\'ees d\'efinissent une m\'etrique
sur~$L$, dite \emph{lisse} (ou \emph{alg\'ebrique}), de la fa\c{c}on
suivante: pour tout ouvert affine~$U$ de~$\mathscr X$ et toute
trivialisation~$\eps$ de~$\mathscr L$ sur~$U$, on a $\norm{\eps}=1$ sur
$\sp^{-1}(U_s)$. Cette d\'efinition est valide car pour toute
fonction~$f_U\in \mathscr{O}_{\mathscr{X}}(U)^{\times}$, la fonction
$\abs{f_{U}}$ sur $\sp_{\mathscr{X}}^{-1}(U_{s})$ est identiquement
\'egale \`a $1$. Nous dirons que $\bar L$ est \emph{semi-positif} si
sa m\'etrique est limite uniforme de m\'etriques lisses donn\'ees
par des mod\`eles~$(\mathscr X,\mathscr L)$ tels que $(c_1(\mathscr
L)|C)\geq 0$ pour toute courbe irr\'eductible $C$ de la fibre
sp\'eciale de~$\mathscr X$.

\medskip

On dit qu'un fibr\'e en droites m\'etris\'e~$\bar L$ est
\emph{admissible} (Zhang dit \emph{int\'egrable}
dans~\cite{zhang95b}) s'il est quotient de deux fibr\'es en droites
munis de m\'etriques semi-positives. 

Toute fonction continue $\phi$ sur $X_v$ d\'efinit une m\'etrique
sur le fibr\'e trivial $\mathscr O_X$, telle que
$\log\norm{1}^{-1}=\phi$. On dira que $\phi$ est lisse, resp.\
admissible, si ce fibr\'e m\'etris\'e $\mathscr O_X(\phi)$ est
lisse, resp.\ admissible.

\subsection*{Fonctions de Green}
Soit $D$ un diviseur de Cartier effectif sur~$X_v$; il lui
correspond un fibr\'e en droites $\mathscr O(D)$ muni d'une section
globale~$s_D$ de diviseur~$D$. On appellera \emph{fonction de Green}
\emph{continue}, resp.\ \emph{lisse}, resp.\ \emph{semi-positive},
resp.\ \emph{admissible} pour~$X_v$ une fonction
$g_D$ sur $X_v\setminus D$ telle que $g_D=\log\norm{s_D}^{-1}$, o\`u
$\norm{\cdot}$ est une m\'etrique continue, resp.\ lisse, resp.\
semi-positive, resp.\ admissible sur~$\mathscr O_X(D)$. On \'etend
par lin\'earit\'e cette terminologie aux diviseurs qui ne sont pas
n\'ecessairement effectifs.

En particulier, \`a tout couple $(\mathscr{X}, \mathscr{D})$,
constitu\'e d'un mod\`ele $\mathscr{X}$ de $X$ sur~$R$
et d'un diviseur $\mathscr{D}$ sur $\mathscr{X}$, est associ\'ee
une fonction de Green lisse $g_{\mathscr{D}}$ pour le diviseur $\mathscr{D}_{|X}$
sur $X_{v}$, caract\'eris\'ee de la fa\c{c}on suivante :
quels que soient l'ouvert affine $U$ de $\mathscr{X}$ 
et la section $f_{U} \in \mathscr{O}_{\mathscr{X}}(U)$
telle que $\mathscr{D}_{|U} = \div(f_{|U})$, la restriction
de $g_{\mathscr{D}}$ au domaine affino\"{\i}de $\sp_{\mathscr{X}}^{-1}(U)$ de
$X_{v}$  est la fonction $\log\abs{f_{U}}^{-1}$.

Par d\'efinition, les \emph{fonctions lisses} sur $X_{v}$ sont
pr\'ecis\'ement les fonctions de Green provenant de couples
$(\mathscr{X}, \mathscr{D})$ tels que $\mathscr{D}_{|X}=0$, \emph{i.e.}
tels que le support de $\mathscr{D}$ soit contenu dans la fibre
sp\'eciale $\mathscr{X}_{s}$ de $\mathscr{X}$.

Il est manifeste qu'un diviseur $\mathscr{D}$ \emph{effectif} donne
lieu \`a une fonction de Green \emph{positive} ; lorsque $\mathscr{X}$ est \emph{int\'egralement ferm\'e} dans sa fibre g\'en\'erique $X$, l'assertion
r\'eciproque est vraie d\`es que $\mathscr{D}_{|X}$ est effectif 
en vertu de la proposition suivante.

\begin{prop}\label{prop.effectif} 
Soit $\mathscr{X}$ un $R$-sch\'ema de
pr\'esentation finie et plat, que l'on suppose int\'egralement ferm\'e dans
sa fibre g\'en\'erique~$X$.

Un diviseur $\mathscr{D}$ sur $\mathscr{X}$
est effectif si et seulement si le diviseur $\mathscr{D}_{|X}$ est
effectif et si la fonction de Green $g_{\mathscr{D}}$ est positive.
\end{prop}

\begin{proof}
Si $\mathscr{D}$ est effectif, il en
est de m\^eme pour $\mathscr{D}_{|X}$ et la fonction de Green
$g_{\mathscr{D}}$ est positive. Supposons r\'eciproquement que
$\mathscr{D}_{|X}$ soit effectif et que la fonction $g_{\mathcal{D}}$
soit positive. \'Etant donn\'e un ouvert affine $U$ de $\mathscr{X}$
tel que $\mathscr{D}_{|U}$ soit principal, il d\'ecoule de nos
hypoth\`eses que $\mathscr{D}_{|U}$ est le diviseur d'un \'el\'ement
r\'egulier $f$ de $\mathcal{O}_{X}(U) = \mathcal{O}_{\mathscr{X}}(U)
\otimes_R K$ tel que $\abs f(x) \leq 1$ pour toute semi-norme multiplicative
born\'ee~$x$ sur~$\mathcal{O}_{X}(U)$. Posant $\mathcal{A} =
\mathcal{O}_{\mathscr{X}}(U)$ et $A = \mathcal{O}_{X}(U) = \mathcal{A}
\otimes_R K$, le morphisme canonique $\mathcal{A} \rightarrow A$
est injectif par platitude et l'anneau $\mathcal{A}$ est int\'egralement
ferm\'e dans $A$. Nous devons montrer que $f$ appartient \`a~$\mathscr A$.

Consid\'erons la sous-$R$-alg\`ebre $\mathcal{A}'$ de $A[f^{-1}]$
engendr\'ee par $f^{-1}$ et $\mathcal{A}$ ; c'est une $R$-alg\`ebre
de type fini et plate telle que $\mathcal{A}' \otimes_R K = A[f^{-1}]$.

Supposons par l'absurde que $f$ n'appartient pas \`a $\mathcal{A}'$.
Alors, $f^{-1}$ est un \'el\'ement non inversible de $\mathcal{A}'$ et
il existe donc un id\'eal maximal $\mathfrak{m}'$ de $\mathcal{A}'$
contenant~$f^{-1}$. Comme $f^{-1}$ est inversible dans
$\mathcal{A}' \otimes_R K$, $\mathfrak{m}'$ contient n\'ecessairement
$\mathfrak{m}\mathcal{A}'$, o\`u $\mathfrak{m}$ d\'esigne l'id\'eal
maximal de $R$, et $\mathfrak{m}'$ d\'efinit donc un point ferm\'e
de la fibre sp\'eciale $U_s$ de $U$.

L'application de r\'eduction $\sp_{\mathcal{X}} :
\sp_{\mathscr{X}}^{-1}(U_s) \rightarrow U_s$ \'etant surjective
(lemme~\ref{lemm.int.cls}), il existe une semi-norme
multiplicative et born\'ee~$x$ sur $\mathcal{O}_{X}(U)$ telle que
$\sp_{\mathscr{X}}(x) = \mathfrak{m}'$ et l'appartenance de $f^{-1}$
\`a $\mathfrak{m}'$ implique alors l'in\'egalit\'e stricte $\abs{f^{-1}}(x)
< 1$. Comme $\abs f(x) \leq 1$ par hypoth\`ese, nous obtenons une
contradiction et en concluons que $f$ appartient \`a~$\mathscr A'$.

Par cons\'equent, $f$ est entier sur~$\mathscr A$,
donc appartient \`a~$\mathscr A$ puisque $\mathscr A$ 
est int\'egralement ferm\'e dans~$A$ et que $f\in A$.
Il en r\'esulte que le diviseur $\mathscr{D}_{|U}
= \div(f)$ est effectif.
\end{proof}

\subsection*{Mesures et hauteurs}
Soit $\bar L_0$, \dots, $\bar L_d$ des fibr\'es en droites
admissibles sur~$X$. Pour tout entier $i$ tel que \mbox{$0\leq i \leq d$},
soit $s_i$ une section de~$\bar L_i$; supposons que les diviseurs de
ces sections se coupent proprement. Pour toute sous-vari\'et\'e~$Z$
de~$X$, de dimension~$e$, et toute place~$v$ de~$F$, on a d\'efini
dans~\cite{chambert-loir2006} un nombre r\'eel
$(\hdiv(s_0)\dots\hdiv(s_e)|Z)_v$, par approximation \`a partir du
cas de m\'etriques lisses, o\`u il co\"{\i}ncide alors avec
la hauteur locale construite par Gubler dans~\cite{gubler1998}. En
outre, ces nombres r\'eels sont des hauteurs locales pour les
hauteurs d\'efinies par Zhang dans~\cite{zhang95b} et associ\'ees
\`a des fibr\'es en droites munis de m\'etriques ad\'eliques
admissibles, au sens o\`u l'on a
\[ \sum_v (\hdiv(s_0)\dots\hdiv(s_e)|Z)_v = (\hc_1(\bar L_0)\dots\hc_1(\bar L_e)|Z), \]
$v$ parcourant l'ensemble des places de~$F$.

La mesure $c_1(\bar L_1)\dots c_1(\bar L_e)\delta_Z$ sur~$X_v$
d\'efinie dans~\cite{chambert-loir2006}
est support\'ee par~$Z_v$.
Si les m\'etriques des fibr\'es~$\bar L_i$ sont lisses,
donn\'ees par des mod\`eles~$\mathscr L_i$
sur un mod\`ele entier~$\mathscr X$ int\'egralement ferm\'e dans~$X$,
et que l'adh\'erence de Zariski~$\mathscr Z$ de~$Z$ dans~$\mathscr X$
est \'egalement int\'egralement ferm\'ee dans~$Z$, on a
\[ c_1(\bar L_1)\dots c_1(\bar L_e)\delta_{Z} \\
   = \sum_{\xi\in \mathscr Z_s^{(0)}} (c_1(\mathscr L_1)\dots c_1(\mathscr L_e)|     [\overline{\xi}])   \delta_{\sp_Z^{-1}(\xi)} .\]
(D'apr\`es le lemme~\ref{lemm.int.cls}, 
un point g\'en\'erique~$\xi$
de~$\mathscr Z_s$ a un unique ant\'ec\'edent dans~$Z_v$ par l'application de
r\'eduction $\sp_Z\colon Z_v\ra \mathscr Z_s$,
car $\mathscr Z$ est int\'egralement ferm\'e dans~$Z$; 
$\delta_{sp_Z^{-1}(\xi)}$ d\'esigne ainsi la mesure de Dirac 
en ce point de~$Z_v$.)
Lorsque $L_0$ est le
fibr\'e en droites  trivial, $\log\norm{s_0}^{-1}$ est une fonction
continue sur~$X_v$ et l'on a la formule
\[ (\hdiv(s_0)\dots\hdiv(s_e)|Z)_v = \int_{X_v} \log\norm{s_0}_v^{-1}
  c_1(\bar L_1)\dots c_1(\bar L_e)\delta_Z. \]
Par la d\'efinition m\^eme des hauteurs locales, on a, si $\div(s_0)=D$,
\begin{multline*} (\hdiv(s_0)\dots\hdiv(s_e)|Z)_v \\
    =
 (\hdiv(s_1)\dots\hdiv(s_e)|Z\cap D)_v
 + \int_{X_v} \log \norm{s_0}_v^{-1}  c_1(\bar L_1)\dots c_1(\bar L_e)\delta_Z,
\end{multline*}
\`a condition que les m\'etriques des fibr\'es~$\bar L_1$, \dots, $\bar L_e$
soient lisses.
La d\'efinition des mesures $c_1(\bar L_1)\dots c_1(\bar L_e)\delta_Z$
est \'etendue par passage \`a la limite
et multilin\'earit\'e au cas de fibr\'es en droites admissibles
arbitraires (cf.~\cite{chambert-loir2006}, \S2.7).

Lorsque $L_i$ est le fibr\'e trivial, muni d'une m\'etrique
admissible telle que $\log\norm{1}^{-1}=\phi_i$,
l'expression $c_1(\bar L_i)$ est not\'ee $\ddc \phi_i$.
Lorsque $L_i=\mathscr O_X(D_i)$ et que sa m\'etrique
est donn\'ee par une fonction de Green admissible~$g_{D_i}$
pour~$D_i$, on \'ecrit $c_1(D_i,g_{D_i})$ pour~$c_1(\bar L_i)$.
Lorsque $v$ est une place archim\'edienne, ces notations sont compatibles
avec  les constructions de mesures en g\'eom\'etrie complexe
dues \`a Bedford--Taylor, Demailly, etc. (voir notamment
\cite{bedford-t82,demailly1985} ainsi que~\cite{maillot2000}).
On a alors $c_1(D,g_D)=\ddc  g_D+\delta_D$ et il s'agit alors
de produits de courants positifs ferm\'es d\'erivant localement
d'un potentiel continu.

La commutativit\'e de l'accouplement de hauteurs locales (\cite{gubler1998},
prop.~9.3)
entra\^{\i}ne que la mesure $c_1(\bar L_1)\dots c_1(\bar L_e)\delta_Z$
est sym\'etrique en les~$\bar L_i$.
Elle entra\^{\i}ne aussi la formule suivante, variante de la formule de Stokes:
\begin{prop}\label{prop.stokes}
Soit $\phi_0$, $\phi_1$ des fonctions admissibles sur~$X_v$,
soit $\bar L_1,\ldots, \bar L_e$  des fibr\'es en droites admissibles
sur~$X$, munis de sections~$s_i$ dont les diviseurs se coupent
proprement.
Alors
\[ \int_{X_v}\phi_0 \ddc \phi_1 c_1(\bar L_2)\dots c_1(\bar L_e)\delta_Z
= \int_{X_v}\phi_1 \ddc \phi_0 c_1 (\bar L_2)\dots c_1(\bar L_e)\delta_Z.\]
\end{prop}

Soit $\bar L$ un fibr\'e en droites sur~$X$ muni d'une m\'etrique
semi-positive.  Lorsque le fibr\'e en droites~$L$ sous-jacent
\`a~$\bar L$ est ample, voire lorsqu'il est gros, on a
$\deg_L(X)>0$. Nous poserons alors $\mu_{\bar L,v}=c_1(\bar
L)_v^d/\deg_L(X)$; c'est une mesure de probabilit\'e sur~$X_v$.

\section{Approximation de fonctions de Green par des fonctions lisses}
\label{sec.approximation}

\begin{theo}\label{theo.approx}
Soit $F$ un corps de nombres et soit $v$ une place de~$F$.
Soit $X$ un sch\'ema projectif int\`egre sur~$F_v$;
soit $D$ un diviseur effectif sur~$X$ et soit $g$ une
fonction de Green pour~$D$, lisse et semi-positive. Il existe une
suite $(g_n)$ de fonctions lisses sur~$X_v$ v\'erifiant les
propri\'et\'es suivantes:
\begin{enumerate}
\item la suite $(g_n)$ converge simplement vers~$g$
et $g_n\leq g$ pour tout~$n$;
\item pour tout~$n$, $g-g_n$ est une fonction de Green pour~$D$,
semi-positive.
\end{enumerate}
\end{theo}

La d\'emonstration consiste \`a v\'erifier que l'on peut
poser $g_n=\min(g,n)$ dans le cas ultram\'etrique,
et une r\'egularisation de cette fonction dans le cas archim\'edien.

\begin{proof}[D\'emonstration dans le cas archim\'edien]
Pour $\eps>0$, soit $\rho_\eps$ la fonction de~$\R$ dans~$\R$
d\'efinie par~$\rho_\eps(x)=\frac12(x+\sqrt{x^2+\eps^2})$.
Elle est $\mathscr C^\infty$ et sa d\'eriv\'ee,
donn\'ee pour $x\neq 0$ par
\[ \rho'_\eps(x)=\frac12 +\frac x{2\sqrt{x^2+\eps^2}} = 
\frac12\left(1+\operatorname{signe(x)}
    \left( 1+ \frac{\eps^2}{x^2}\right)^{-1/2}\right), \]
est strictement croissante et comprise entre~$0$ et~$1$.
En particulier, $\rho_\eps$ est strictement croissante
et strictement convexe.
En outre, pour tout $x\in\R$, on a l'in\'egalit\'e
\[ \abs{\rho_1(x)-\max( x,0)} = \frac12 \left( x+\sqrt{x^2+1}\right)
 - \frac12 \left(x+\abs x\right) 
 \leq \frac1{\max(1,\abs x)}. \]
Comme $\rho_\eps(x)=\eps\rho_1(x/\eps)$, cela entra\^{\i}ne la majoration
\[ \abs{\rho_\eps(x)-\max( x,0)} 
\leq  \frac{\eps  }{ \max(1,\abs{x/\eps})} \leq \eps.\]
Par suite, lorsque $\eps\ra 0$,
$\rho_\eps$ converge uniform\'ement, en d\'ecroissant, sur~$\R$ vers
la fonction~$x\mapsto \max( x,0)$.

Pour $n~\in~\N$, $\eps>0$ et $x\in\R$, posons
\[ \phi_{n,\eps}(x) = x - \rho_\eps(x-n). \]
Lorsque $\eps\ra 0$, $\phi_{n,\eps}(x)$
tend uniform\'ement en~$n$ et~$x$ vers
$x-\max(x-n,0)=x+\min(n-x,0)=\min(n,x)$.
En outre, la fonction $x\mapsto \phi_{n,\eps}$ est strictement
croissante et strictement concave. On peut alors choisir une suite
$(\eps_n)$ de nombres r\'eels qui tende vers~$0$ en d\'ecroissant de
sorte que la suite de fonctions~$(\phi_{n,\eps_n})$ soit une suite
croissante de fonctions $\mathscr C^\infty$ qui converge vers la
fonction $x\mapsto x$, uniform\'ement sur chaque intervalle major\'e. 

Pour $n\in\N$, posons $g_n=\phi_{n,\eps_n}\circ g$; c'est
une fonction lisse sur~$X_v$ et $g_n(x)\leq g(x)$ pour tout
$x\in\R$. De plus, lorsque $n\ra\infty$, $g_n(x)$ tend vers~$g(x)$.
Pour $x\in X_{v}$, on a
\[ g(x)-g_n(x) =g(x)-\phi_{n,\eps_n}(g(x))=\rho_{\eps_n}(g(x)-n). \]
C'est une fonction de Green lisse, positive, pour~$D$. Comme la
fonction $\rho_{\eps_n}$ est convexe et strictement croissante,
$g-g_n$ est semi-positive. Cela termine la d\'emonstration lorsque
$v$ est une place archim\'edienne de~$F$.
\end{proof}

\begin{proof}[D\'emonstration du th\'eor\`eme dans le cas d'une place ultram\'etrique]
Nous consid\'erons maintenant une place non archim\'edienne $v$ de
$F$. Soit $R$ l'anneau des entiers du corps local $F_{v}$ et soit
$\pi$ une uniformisante de $R$.  Pour tout entier naturel~$n$,
posons $g_{n} = \min(g, n \log|\pi|^{-1})$; il s'agit d'une fonction
continue et born\'ee sur~$X_{v}$, et nous allons v\'erifier que
$g-g_{n}$ est une fonction de Green semi-positive et lisse pour~$D$.
La suite croissante~$(g_{n})$ converge simplement vers~$g$.

Soit $\mathscr X$ un mod\`ele entier de~$X$ int\'egralement ferm\'e dans sa fibre g\'en\'erique et soit $\mathscr D$
un diviseur de Cartier sur~$\mathscr X$ qui prolonge~$D$
et qui d\'efinit la fonction de Green lisse~$g$.
Comme $g$ est minor\'ee, il existe un entier $n_0$ 
tel que $g+n_0\log\abs{\pi}^{-1}$ soit positive.
C'est une fonction de Green pour~$D$, 
associ\'ee au diviseur $\mathscr D+n_0[\mathscr X_s]$;
comme $\mathscr X$ est int\'egralement ferm\'e dans $X$, ce diviseur est effectif.
Quitte \`a remplacer $g$ par $g+n_{0}\log|\pi|^{-1}$ et $n$ par
$n_{0} + n$, nous pouvons donc supposer que $\mathscr D$ est effectif  ; 
c'est ce que nous faisons.

Notons $\mathscr{I}$ le faisceau coh\'erent d'id\'eaux
d\'efini par
\[\mathscr{I} = \mathscr{I}_{\mathscr{X}}(D) + \pi^{n}\mathscr{O}_{\mathscr{X}},\]
dont le support est contenu dans la fibre sp\'eciale de
$\mathscr{X}_{s}$. Soit $p \colon \mathscr{X}' \rightarrow
\mathscr{X}$ l'\'eclatement de $\mathscr{I}$ et soit $\mathscr{E}$
le diviseur exceptionnel, d\'efini par le faisceau d'id\'eaux
$p^{-1}\mathscr I\cdot\mathscr O_{\mathscr X'}$.
Le diviseur
\[\mathscr{D}' = p^{*}\mathscr{D} - \mathscr{E}\] prolonge~$D$ sur
le mod\`ele~$\mathscr{X}'$ de $X$ et d\'efinit donc une fonction de
Green~$g'$ pour~$D$ sur~$X_{v}$.

\'Etant donn\'e un ouvert affine~$V$ de~$\mathscr{X}$
tel que $\mathscr{D}_{|V} = \div(f_{V})$, pour $f_{V} \in
\mathscr{O}_{\mathscr{X}}(V)$, le $V$-sch\'ema~$\mathscr{X}'_{|V} =
p^{-1}(V)$ s'identifie naturellement au sous sch\'ema ferm\'e du
$V$-sch\'ema projectif d\'efini par l'\'equation $f_{V}T_{0}-\pi^{n}T_{1}=0$
dans $\P^{1}_{V} = \Proj(\mathscr{O}_{V}[T_{0}, T_{1}])$ ;
$\mathscr D'$ est alors d\'efini par l'annulation de~$T_1$, $\mathscr E$
par celle de~$\pi^n$. 
Il en d\'ecoule que $\mathscr{X}'$ est recouvert par des
ouverts affines~$U$ tels que $p^{*}\mathscr{D}_{|U} = \div(f'_{U})$,
pour $f'_{U} \in \mathscr{O}_{\mathscr{X}'}(U)$, et satisfaisant \`a
l'une des deux conditions suivantes:
\begin{itemize} \item ou bien $f'_{U} \in
\pi^{n}\mathscr{O}_{\mathscr{X}'}(U)$, auquel cas $\mathscr{E}_{|U}
= \div(\pi^{n})$; \item ou bien $\pi^{n} \in
f'_{U}\mathscr{O}_{\mathscr{X}'}(U)$, auquel cas $\mathscr{E}_{|U}
= \div(f'_{U}) = p^{*}\mathscr{D}_{|U}$.  \end{itemize}

Un tel recouvrement ouvert $(U)$ de $\mathscr{X}$ induit un recouvrement
de l'espace analytique $X_{v}$ par les domaines affino\"{\i}des
$\sp^{-1}(U\cap\mathscr X_s)$ 
sur lesquels ou bien 
$g = \log\abs{f'_{U}}^{-1} \geq n\log\abs\pi^{-1}$, ou bien
$g = \log\abs{f'_{U}}^{-1} \leq n\log\abs\pi^{-1}$. 
Par cons\'equent, $g'$
co\"{\i}ncide avec $g-g_{n} = \max(g-n\log\abs\pi^{-1}, 0)$ sur
$X_{v}$ et cette derni\`ere est donc bien une fonction de Green
lisse pour $D$. Il reste \`a v\'erifier que le diviseur
$\mathscr{D}'$ est nef sur la fibre sp\'eciale de $\mathscr{X}'$.

La description locale qui pr\'ec\`ede fournit
imm\'ediatement l'\emph{effectivit\'e} du diviseur $\mathscr{D}'$.
En outre, d\'esignant par $\Omega$ le plus grand ouvert de $\mathscr{X}'$
sur lequel $\pi^{n}$ engendre $p^{-1}\mathscr{I}$, il est clair que
le support de $\mathscr{D}'$ est contenu dans $\Omega$ et que
$\mathscr{E}_{|\Omega} = \div(\pi^{n})$.

Consid\'erons alors un sous-sch\'ema ferm\'e int\`egre~$C$, de
dimension un et contenu dans la fibre sp\'eciale~$\mathscr{X}'_{s}$
de $\mathscr{X}'$. Si $C$ est g\'en\'eriquement disjoint du support
de $\mathscr{D}'$, on peut d\'efinir la
restriction~$\mathscr{D}'_{|C}$ de~$\mathscr D'$ \`a~$C$ ; c'est un
diviseur effectif sur $C$ et le nombre d'intersection
\[(c_{1}(\mathscr{D}') | C) = \deg(\mathscr{D}'_{|C})\] est donc
positif. Si, au contraire, $C$ est contenu dans le support de
$\mathscr{D}'$, $\Omega$ est alors un voisinage de $C$ sur lequel
les diviseurs $\mathscr{D}'$ et $p^{*}\mathscr{D}$ sont
\'equivalents et, en vertu de la formule de projection, le nombre
d'intersection \[(c_{1}(\mathscr{D}')|C) =
(c_{1}(p^{*}\mathscr{D})|C) = (c_{1}(\mathscr{D}) | p_{*}C)\] est
positif puisque, par hypoth\`ese, le diviseur $\mathscr{D}$ est nef
sur $\mathscr{X}_{s}$. Nous venons de v\'erifier que le diviseur
$\mathscr{D}'$ est nef sur la fibre sp\'eciale de $\mathscr{X}'$ et
ceci ach\`eve la d\'emonstration du th\'eor\`eme.
\end{proof}

\begin{coro}\label{coro.approx}
Reprenons les notations du th\'eor\`eme et consid\'erons
des fibr\'es en droites admissibles $\bar L_1$,\dots, $\bar L_{d-1}$
sur~$X$.
La suite $(\ddc g_nc_1(\bar L_1)\dots c_1(\bar L_{d-1}))$
de mesures sur~$X_v$ est relativement compacte dans l'espace
des mesures born\'ees sur~$X_v$, pour la topologie \'etroite.
En particulier, pour montrer qu'elle converge vers une mesure~$\mu$,
il suffit de montrer que pour toute fonction lisse~$\phi$,
\[ \int_{X_v} \phi \ddc g_nc_1(\bar L_1)\dots c_1(\bar L_{d-1})
 \xrightarrow[n\ra\infty]{} \mu(\phi). \]
\end{coro}
\begin{proof}
On peut supposer que les fibr\'es~$\bar L_i$ sont semi-positifs.
Pour tout $n$, l'\'egalit\'e
\[  \ddc g_nc_1(\bar L_1)\dots c_1(\bar L_{d-1})
 = c_1(D,g) c_1(\bar L_1)\dots c_1(\bar L_{d-1})
   - \ddc (g-g_n) c_1(\bar L_1)\dots c_1(\bar L_{d-1}) \]
exprime le premier membre comme diff\'erence de deux mesures
positives dont la masse totale est \'egale \`a
\[  (c_1(L_1)\dots c_1(L_{d-1})|D). \]
La premi\`ere assertion d\'ecoule donc de la compacit\'e faible
de l'ensemble des mesures de probabilit\'e sur l'espace compact~$X_v$.
La seconde assertion du corollaire est alors cons\'equence
de ce que l'ensemble des fonctions lisses sur~$X_v$ est dense
dans l'ensemble des fonctions continues (\cf\cite{gubler1998}, th.~7.12).
\end{proof}

\section{Int\'egration de fonctions de Green}
\label{sec.mahler}

Nous allons en fait
d\'emontrer une version locale du th\'eor\`eme~\ref{theo.mahler}.

\begin{theo}\label{ineg.mahler}
Fixons une place~$v$ de~$F$.
Soit $D$ un diviseur de Cartier sur~$X$ et soit $g_D$ une fonction
de Green admissible pour $D$; on d\'esigne par $\bar L_0$ le fibr\'e
en droites $\mathscr O_X(D)$ muni de la m\'etrique admissible pour
laquelle la norme de la section canonique~$s_0$ (de diviseur~$D$)
v\'erifie $\log\norm{s_0}_v^{-1}=g_D$. Soit, pour $1\leq i\leq d$, un
fibr\'e en droites admissible $\bar L_i$ sur~$X_v$ muni d'une
section globale~$s_i$. On suppose que les diviseurs
de~$s_0,\dots,s_d$ se coupent proprement.

Alors la fonction~$g_D$ est int\'egrable
pour la mesure (sign\'ee) $c_1(\bar L_1)\dots c_1(\bar L_d)$ et l'on a
\[ \int_{X_v} g_D  c_1(\bar L_1)\dots c_1(\bar L_d)
   = (\hdiv(s_0)\dots\hdiv(s_d)|X)_v - (\hdiv(s_1)\dots\hdiv(s_d)|D)_v. \]
\end{theo}
Faisant varier~$v$ parmi les places de~$F$ et additionnant
les formules obtenues, on en d\'eduit
ais\'ement le th\'eor\`eme~\ref{theo.mahler}, au moins
lorsque $\bar L$ est tr\`es ample (de sorte qu'il poss\`ede
des sections qui se coupent proprement sur~$D$);
le cas g\'en\'eral s'ensuit par multilin\'earit\'e.

Avant de passer \`a la d\'emonstration du th\'eor\`eme~\ref{ineg.mahler},
donnons-en un corollaire,
dont la version aux places archim\'ediennes est bien connue.
\begin{coro}\label{coro.ineg.mahler}
Si $Z$ est un ferm\'e irr\'eductible de~$X$, de dimension~$e$, et
$v$ est une place de~$F$, la mesure $c_1(\bar L_1)\dots c_1(\bar
L_e)\delta_Z$ sur~$X_v$ ne charge pas les sous-ensembles
alg\'ebriques stricts de~$Z$.
\end{coro}
\begin{proof}
Il suffit de traiter le cas $Z=X$. Tout sous-ensemble alg\'ebrique
strict de~$X$ est contenu dans un diviseur de Cartier effectif~$D$;
soit $g_D$ une fonction de Green admissible pour~$D$.
La fonction $g_D$ est int\'egrable contre
la mesure $c_1(\bar L_1)\dots c_1(\bar L_d)$ et vaut~$+\infty$ sur~$D$.
N\'ecessairement, $D$ est de mesure nulle pour $c_1(\bar L_1)\dots c_1(\bar L_d)$,
et $Z$ aussi.
\end{proof}

D\'emontrons maintenant le th\'eor\`eme~\ref{ineg.mahler}.
(Pour all\'eger les notations, nous omettrons le plus souvent les indices~$v$.)

Lorsque $v$ est une place archim\'edienne,
l'int\'egrabilit\'e et l'\'egalit\'e du th\'eor\`eme~\ref{ineg.mahler}
 sont de nature locale ; elles d\'ecoulent alors
des propri\'et\'es de continuit\'e de l'op\'erateur de
Monge-Amp\`ere complexe, voir par exemple le th\'eor\`eme~2.2
de~\cite{demailly1985}. La d\'emonstration que nous allons donner
ici -- qui ne n\'ecessite pas de distinguer les places
archim\'ediennes des places ultram\'etriques -- est essentiellement
formelle une fois que le th\'eor\`eme d'approximation est acquis.

\medskip

\emph{R\'eduction au cas o\`u les~$\bar L_i$ sont semi-positifs.}
--- Par multilin\'earit\'e, il est loisible de supposer que les
fibr\'es en droites m\'etris\'es~$\bar L_i$ sont semi-positifs. La
mesure $c_1(\bar L_1)\dots c_1(\bar L_d)$ est alors positive. Fixons
aussi, pour tout~$i\in\{1,\dots,d\}$, une m\'etrique lisse
semi-positive sur~$L_i$, correspondant \`a un fibr\'e en droites
m\'etris\'e $\bar M_i$. Soit $t_i$ la section de~$M_i$ correspondant
\`a~$s_i$ et soit $\phi_i$ la fonction admissible $\log
\norm{t_i}^{-1} - \log\norm{s_i}^{-1}$, de sorte que $\bar L_i=\bar
M_i\otimes\mathscr O(\phi_i)$.

\medskip

\emph{R\'eduction au cas o\`u $D$ est ample et $g_D$ lisse,
semi-positive.} --- Par lin\'earit\'e, il suffit de traiter le cas
o\`u $D$ est un diviseur ample sur~$X$. Soit alors $g_D^0$ une
fonction de Green lisse, semi-positive pour~$D$, de sorte que
$\phi=g_D-g_D^0$ est une fonction admissible sur~$X_v$. Si le
th\'eor\`eme est vrai pour~$g_D^0$, alors $g_D=\phi+g_D^0$ est
int\'egrable contre $c_1(\bar L_1)\dots c_1(\bar L_d)$ et
\begin{align*}
 \int_{X_v} g_D c_1(\bar  L_1)\dots c_1(\bar L_d)
& =  \int_{X_v} g_D^0 c_1(\bar L_1)\dots c_1(\bar L_d)
+ \int_{X_v} \phi c_1(\bar L_1)\dots c_1(\bar L_d) \\
&= (\hdiv(s_0)^0\hdiv(s_1)\dots\hdiv(s_d)|X)
  - (\hdiv (s_1)\dots\hdiv(s_d)|D)  \\
& \hskip 6cm + \int_{X_v} \phi c_1(\bar L_1)\dots c_1(\bar L_d),
\end{align*} o\`u $\hdiv(s_0)^0$ correspond \`a la section
canonique~$s_0$ du fibr\'e en droites~$\mathscr O_X(D)$ muni de la
m\'etrique lisse donn\'ee par~$g_D^0$. Par d\'efinition de la
mesure~$c_1(\bar L_1)\dots c_1(\bar L_d)$, on a
\begin{multline*} (\hdiv(s_0)^0\hdiv(s_1)\dots\hdiv(s_d)|X) \\
= (\hdiv(s_0)\hdiv(s_1)\dots\hdiv(s_d)|X)
 - \int_{X_v} \phi c_1(\bar L_1)\dots c_1(\bar L_d) , \end{multline*}
si bien que
\[ \int_{X_v} g_D c_1(\bar  L_1)\dots c_1(\bar L_d)
 = (\hdiv(s_0)\hdiv(s_1)\dots\hdiv(s_d)|X)
      - (\hdiv (s_1)\dots\hdiv(s_d)|D), \]
ce qu'il fallait d\'emontrer.

On suppose donc dans la suite que $D$ est ample et que $g_D$ est une
fonction de Green lisse, semi-positive pour~$D$.

\emph{Montrons par r\'ecurrence sur $k\in\{0,\dots,d\}$
que l'\'enonc\'e vaut lorsque~$\bar L_i$ 
est muni d'une m\'etrique lisse pour $i>k$.} 

L'assertion est vraie pour $k=0$ (on est dans le cas lisse).
Supposons-la vraie pour $k-1$ et d\'emontrons-la pour~$k$.
Il suffit de d\'emontrer l'\'egalit\'e lorsque $\bar L_i=\bar M_i$ pour $i>k$.
Comme $\bar L_k =\bar M_k\otimes\mathscr O(\phi_k)$, on a l'\'egalit\'e de mesures sur~$X_v$:
\begin{multline*}
c_1(\bar L_1)\dots c_1(\bar L_k)c_1(\bar M_{k+1})\dots c_1(\bar M_d) \\
= c_1(\bar L_1)\dots c_1(\bar L_{k-1})c_1(\bar M_{k})\dots c_1(\bar M_d)\\
+ \ddc \phi_k c_1(\bar L_1)\dots c_1(\bar L_{k-1})c_1(\bar M_{k+1})\dots c_1(\bar M_d)
\end{multline*}
et les \'egalit\'es de hauteurs locales:
\begin{multline*}
(\hdiv(s_0)\hdiv(s_1)\dots\hdiv(s_{k})\hdiv(t_{k+1})\dots\hdiv(t_d)|X)\\
= (\hdiv(s_0)\hdiv(s_1)\dots\hdiv(s_{k-1})\hdiv(t_k)\dots\hdiv(t_d)|X)\\
+\int_{X_v} \phi_k c_1(D,g_D) c_1(\bar L_1)\dots c_1(\bar L_{k-1})
c_1(\bar M_{k+1})\dots c_1(\bar M_d)
\end{multline*}
et
\begin{multline*}
(\hdiv(s_1)\dots\hdiv(s_{k})\hdiv(t_{k+1})\dots\hdiv(t_d)|D)\\
= (\hdiv(s_1)\dots\hdiv(s_{k-1})\hdiv(t_k)\dots\hdiv(t_d)|X)\\
+\int_{X_v} \phi_k c_1(\bar L_1)\dots c_1(\bar L_{k-1}) c_1(\bar
M_{k+1})\dots c_1(\bar M_d)\delta_D.
\end{multline*}
Comme $\bar M_k,\dots,\bar M_d$ sont lisses, il r\'esulte de l'hypoth\`ese de
r\'ecurrence que l'on a
\begin{multline*}
 \int_ {X_v} g_D c_1(\bar L_1)\dots c_1(\bar L_{k-1})c_1(\bar M_k)\dots
c_1(\bar M_d) \\
 = (\hdiv(s_0)\hdiv(s_1)\dots \hdiv (s_{k-1})\hdiv(t_k)\dots\hdiv(t_d)|X)\\
- (\hdiv(s_1)\dots \hdiv (s_{k-1})\hdiv(t_k)\dots\hdiv(t_d)|D).
\end{multline*}
Par suite,
\begin{align*}
\int_{X_v} g_D c_1(\bar L_1)\dots c_1(\bar L_k)c_1(\bar M_{k+1})\dots c_1(\bar M_d) \hskip -5cm \\
& = \int_{X_v} g_D c_1(\bar L_1)\dots c_1(\bar L_{k-1})c_1(\bar M_k)c_1(\bar M_{k+1})\dots c_1(\bar M_d) \\
& \qquad + \int_{X_v} g_D \ddc \phi_k c_1(\bar L_1)\dots c_1(\bar L_{k-1})c_1(\bar M_{k+1})\dots c_1(\bar M_d) \\
& = (\hdiv(s_0)\hdiv(s_1)\dots\hdiv(s_{k-1}) \hdiv(t_k)\dots\hdiv(t_d)|X)\\
& \qquad - (\hdiv(s_1)\dots\hdiv(s_{k-1}) \hdiv(t_k)\dots\hdiv(t_d)|D) \\
& \qquad + \int_{X_v} g_D \ddc \phi_k c_1(\bar L_1)\dots c_1(\bar L_{k-1})c_1(\bar M_{k+1})\dots c_1(\bar M_d) \\
& =  (\hdiv(s_0)\hdiv(s_1)\dots\hdiv(s_{k-1}) \hdiv(t_k)\dots\hdiv(t_d)|X)\\*
& \qquad + \int_{X_v} \phi_k c_1(D,g_D) c_1(\bar L_1)\dots c_1(\bar L_{k-1})
c_1(\bar M_{k+1})\dots c_1(\bar M_d) \\*
& \qquad - (\hdiv(s_1)\dots\hdiv(s_{k-1}) \hdiv(t_k)\dots\hdiv(t_d)|D) \\*
& \qquad - \int_{X_v}\phi_k c_1(\bar L_1)\dots c_1(\bar L_{k-1})
               c_1(\bar M_{k+1})\dots c_1(\bar M_d)\delta_D \\*
& \qquad + \int_{X_v} g_D \ddc \phi_k c_1(\bar L_1)\dots c_1(\bar L_{k-1}) c_1(\bar M_{k+1})\dots c_1(\bar M_d).
\end{align*}
Il suffit donc de d\'emontrer l'\'egalit\'e
\begin{multline*}
\int_{X_v} \phi_k c_1(D,g_D)c_1(\bar L_1)\dots c_1(\bar L_{k-1})c_1(\bar M_{k+1})\dots c_1(\bar M_d) \\
{} = \int_{X_v} \phi_k c_1(\bar L_1)\dots c_1(\bar L_{k-1}) c_1(\bar M_{k+1})\dots c_1(\bar M_d) \delta_D \\
{} - \int_{X_v} g_D \ddc \phi_k c_1(\bar L_1)\dots c_1(\bar L_{k-1})c_1(\bar M_{k+1})\dots c_1(\bar M_d).
\end{multline*}

\medskip
\emph{Application du th\'eor\`eme d'approximation.} --- Soit $(g_n)$
une suite de fonctions lisses sur~$X_v$ telle que $g_n\leq g_D$,
$g_n$ converge simplement vers~$g_D$ et telle que $h_n=g_D-g_n$ soit
une fonction de Green semi-positive pour~$D$. Une telle suite existe
d'apr\`es le th\'eor\`eme~\ref{theo.approx}.

Alors,
\begin{align*}
\int_{X_v} g_D \ddc\phi_k c_1(\bar L_1)\dots c_1(\bar L_{k-1})c_1(\bar M_{k+1})\dots c_1(\bar M_d) \hskip-5cm \\
&= \lim_{n\ra\infty} \int_{X_v} g_n \ddc \phi_k c_1(\bar L_1)\dots c_1(\bar L_{k-1})c_1(\bar M_{k+1})\dots c_1(\bar M_d) \\
& = \lim_{n\ra\infty} \int_{X_v} \phi_k \ddc g_n c_1(\bar L_1)\dots c_1(\bar L_{k-1})c_1(\bar M_{k+1})\dots c_1(\bar M_d)
\end{align*}
d'apr\`es la prop.~\ref{prop.stokes}.
La fonction~$\phi_k$ est continue, car admissible; il suffit
donc de montrer la convergence des mesures $\ddc g_nc_1(\bar L_1)\dots c_1(\bar L_{k-1})c_1(\bar M_{k+1})\dots c_1(\bar M_d)$ vers la mesure
\begin{multline*}
 c_1(D,g_D) c_1(\bar L_1)\dots c_1(\bar L_{k-1})c_1(\bar M_{k+1})
 \dots c_1(\bar M_d) \\ 
{} - c_1(\bar L_1)\dots c_1(\bar L_{k-1})c_1(\bar M_{k+1})\dots c_1(\bar M_d)\delta_D,
\end{multline*}
lorsque $n\ra\infty$. D'apr\`es le corollaire~\ref{coro.approx}
du th\'eor\`eme d'approximation,
il suffit de montrer la convergence apr\`es int\'egration
contre une fonction lisse, ce qui revient \`a montrer
l'\'egalit\'e voulue sous l'hypoth\`ese suppl\'ementaire que $\phi_k$
est une fonction lisse, laquelle vaut par l'hypoth\`ese de r\'ecurrence.

Cela termine la d\'emonstration du th\'eor\`eme~\ref{ineg.mahler}.

Le th\'eor\`eme~\ref{theo.mahler} en est une 
cons\'equence \'evidente si $L$ poss\`ede
des sections globales $s_1,\ldots,s_d$
telles que $D$ et les diviseurs des~$s_i$ se coupent proprement:
il suffit alors de sommer les \'egalit\'es fournies
par le th\'eor\`eme~\ref{ineg.mahler} \`a chaque place du corps~$F$.
Le cas g\'en\'eral en r\'esulte alors par multilin\'earit\'e.

\section{Volumes en g\'eom\'etrie d'Arakelov}
\label{sec.volumes}

Ce paragraphe est consacr\'e \`a la reformulation de r\'esultats
de Zhang~\cite{zhang95b} et Yuan~\cite{yuan2006}.

Soit $X$ un sch\'ema projectif purement de dimension $d$ sur un
corps de nombres~$F$ et soit $\bar L$ un fibr\'e en droites muni
d'une m\'etrique ad\'elique admissible. Le volume de~$L$ est
d\'efini par la formule
\[ \vol(L) = \limsup_{k\ra\infty} \frac{\dim_F H^0(X,L^{\otimes k})}{k^{d}/d!} \]
et mesure la croissance du nombre de sections globales des puissances
tensorielles de~$L$. Si $L$ est ample, le th\'eor\`eme
de Hilbert-Samuel implique que $\vol(L)=\deg_L(X)$;
plus g\'en\'eralement, cette \'egalit\'e vaut
si $L$ est num\'eriquement effectif.
On dit que $L$ est gros si son volume est strictement positif.
On d\'eduit facilement d'un th\'eor\`eme de Fujita~\cite{fujita1994}
d'existence de {\og d\'ecompositions de Zariski approch\'ees\fg}
que la limite sup\'erieure est en fait une limite
(\cite{lazarsfeld2004}, 10.3.11).

Pour toute place~$v$ de~$F$, munissons l'espace vectoriel $H^0(X,L)$
des sections globales de~$X$ de la norme $v$-adique
$\norm{s}_v=\sup_{x\in X(\bar F_v)}\norm{s(x)}_v$, pour $s\in H^0(X,
L)$. Soit $\mathbf A_F$ l'anneau des ad\`eles de~$F$ et soit $\mu$
une mesure de Haar sur le groupe ab\'elien localement compact
$H^0(X, L)\otimes\mathbf A_F$. On pose alors
(voir~\cite{zhang95b,bombieri-v1983}):
\[ \chi_{\sup}(\bar L) = - \log \frac{\mu ( H^0(X, L)\otimes\mathbf A_F / H^0(X, L) )}{\mu (\prod_v B_v)}, \]
o\`u $B_v$ d\'esigne la boule unit\'e de $H^0(X, L)\otimes F_v$ pour
la norme~$v$-adique introduite. C'est un analogue en g\'eom\'etrie
d'Arakelov de la dimension de l'espace des sections globales de~$
L$. Poussant l'analogie, on d\'efinit le \emph{volume arithm\'etique} de~$\bar L$
par l'expression
\[ \hvol(\bar L) =
\limsup_{k\ra\infty} \frac{\chi_{\sup} (\bar L^{\otimes k})}{k^{d+1}/(d+1)!}. \] 
D'apr\`es un th\'eor\`eme de Rumely, Lau et Varley~\cite{rumely-l-v2000}, 
c'est en fait une limite, lorsque $L$ est ample. 
(Pour une meilleure analogie, le volume arithm\'etique devrait \^etre d\'efini
en termes de sections globales de normes au plus~$1$ \`a toutes les places,
comme dans~\cite{moriwaki2006}, 
mais la quantit\'e introduite ici est plus maniable.)
Si $\bar L$ est fibr\'e en droites ample muni d'une m\'etrique ad\'elique 
semi-positive, la variante de~\cite{zhang95b}
du th\'eor\`eme de Hilbert-Samuel arithm\'etique
(d\'eduite dans \emph{loc. cit.} de la version de~\cite{gillet-s88,abbes-b95})
affirme que l'on a $ \hvol(\bar L)= h_{\bar L}(X) $.

Soit $v$ une place de~$F$, soit $\phi$ une fonction continue sur~$X_v$
et soit $\bar M=\bar L(\phi)$. Si $s$ est une section
globale de~$L$, on a $\norm{s}^{\bar M}_v(x)=\norm{s}^{\bar L}(x) e^{-\phi(x)}$;
par suite,
\[ \chi_{\sup}(\bar L)+\dim H^0(X,L) \inf\phi
 \leq \chi_{\sup}(\bar M)
 \leq \chi_{\sup}(\bar L) + \dim H^0(X,L) \sup\phi .\]
Par suite, 
\begin{equation*}
   \hvol(\bar L) + (d+1) \vol(L) \inf\phi \leq \hvol(\bar M)
 \leq \hvol(\bar L)+(d+1)\vol(L)\sup\phi. \end{equation*}
En particulier, le volume  arithm\'etique d\'epend contin\^ument de la m\'etrique ad\'elique sur~$\bar L$. 
(Le r\'esultat analogue pour la d\'efinition alternative
du volume arithm\'etique en
termes de sections globales est le th\'eor\`eme principal de~\cite{moriwaki2006}.)

\medskip

Nous dirons qu'un fibr\'e en droites ample~$L$ sur~$X$,
muni d'une m\'etrique ad\'elique semi-positive est
\emph{arithm\'etiquement num\'eriquement effectif},
resp.~\emph{arithm\'etiquement ample},
si l'on a $h_{\bar L}(Y)\geq 0$
pour toute sous-vari\'et\'e~$Y$ de~$X$,
resp. si l'in\'egalit\'e stricte est satisfaite.
(D'apr\`es Zhang~\cite{zhang95b}, il suffirait  d'imposer
cette condition pour les sous-vari\'et\'es de dimension~$0$.)

Si $\bar L$ est semi-positif (resp.~arithm\'etiquement
num\'eriquement effectif), alors $\bar L(c)$ est arithm\'etiquement
ample pour tout nombre r\'eel~$c$ assez grand
(resp. pour tout nombre r\'eel~$c>0$).
Si $\bar L$ est arithm\'etiquement ample,
pour tout entier~$n$ assez grand,
$H^0(X,L^{\otimes n})$ poss\`ede une base form\'ee de sections globales
dont les normes sup. \`a toutes les places
sont major\'ees par~$1$. 
Inversement, lorsque $\bar L$ est semi-positif,
cette condition  entra\^{\i}ne que toute sous-vari\'et\'e de~$X$
est de hauteur positive (voir par exemple~\cite{zhang95b}).

\medskip

Le volume de~$\bar L$ contr\^ole la hauteur des sous-vari\'et\'es 
de~$X$ au sens suivant:
\begin{lemm}\label{lemm.inegalite}
Supposons que $L$ soit gros.
Le minimum essentiel $e(\bar L)$ des hauteurs des points de~$X(\bar F)$
relativement \`a~$\bar L$ v\'erifie l'in\'egalit\'e
\[ e(\bar L) = \sup_{\substack{\emptyset\neq U\subset X \\ \text{\upshape $U$ ouvert}}}
\inf_{x\in U(\bar F)} h_{\bar L}(x)   \geq \frac{\hvol(\bar L)}{(d+1)\vol(L)}. \]
\end{lemm}
\begin{proof}
Soit $\eps$ un nombre r\'eel compris entre~$0$ et~$1$
et soit $n$ un entier tel que
\[ (1-\eps)\vol(L)\frac{k^d}{d!} \leq \dim H^0(X,L^{\otimes k})
 \leq (1+\eps)\vol(L) \frac{k^d}{d!} \]
pour tout entier $k\geq n$ (en particulier, $H^0(X,L^{\otimes k})\neq 0)$).
Par d\'efinition de~$\hvol(\bar L)$, il existe alors un entier~$k\geq n$
tel que
\[ \chi_{\sup} (\bar L^{\otimes
k}) \geq \big(\hvol(\bar L)-\eps\big)
     \frac{k^{d+1}}{(d+1)!}.\]
Notons $D=\dim_F H^0(X,L^{\otimes k})$. 
D'apr\`es le th\'eor\`eme de Minkowski, ou plut\^ot
sa variante ad\'elique due \`a
Bombieri--Vaaler~\cite{bombieri-v1983}, il existe une section non
nulle $\sigma\in H^0(X,L^{\otimes k})$ tel que $\norm{\sigma}_v \leq
1$ pour toute place~$v$ ultram\'etrique et $\norm{\sigma}_v\leq c$
pour toute place~$v$ archim\'edienne, d\`es que le nombre r\'eel~$c$
v\'erifie l'in\'egalit\'e $c^{D}\exp(\chi_{\sup}(\bar L^{\otimes k}))>2^{D}$.
Pour de tels~$c$ et~$\sigma$, et pour tout $x\in X(\bar F)$ tel que
$\sigma(x)\neq 0$, on a
\[ h_{\bar L}(x) 
= - \frac1{[F(x):\Q] k} \sum_{v\in M_F} \sum_{\iota\colon F(x)\hra\bar F_v}
                \log \norm{\sigma}_v(\iota(x))
\geq  -\frac1k \log c ,\] d'o\`u $e(\bar L)\geq - \frac{1}{k}\log c$.
Il en r\'esulte que
\[
 e(\bar L)  \geq -\frac1k \log 2 + \frac1{Dk} \chi_{\sup}(\bar L^{\otimes k}) 
 \geq -\frac1k\log2+\frac{(\hvol(\bar L)-\eps) k^d}{(d+1)! D} . \]
Compte-tenu de l'encadrement pour~$D$, on a donc
\[ e(\bar L) \geq -\frac1k\log 2+ \min\left(\frac{\hvol(\bar L)-\eps}{(d+1)\vol(L)(1-\eps)}, \frac{\hvol(\bar L)-\eps}{(d+1)\vol(L)(1+\eps)}\right).
\]
Lorsque $\eps$ tend vers~$ 0$, $k$ tend vers l'infini et
on obtient l'in\'egalit\'e 
\[ e(\bar L) \geq \frac{\hvol(\bar L)}{(d+1)\vol(L)}, \]
ce qu'il fallait d\'emontrer.
\end{proof}

La fin de ce paragraphe explique le comportement du volume
et du volume arithm\'etique lorsqu'on ajoute \`a un fibr\'e
muni d'une m\'etrique ad\'elique semi-positive un petit multiple
d'un fibr\'e en droites admissible.

Nous aurons \`a utiliser le r\'esultat suivant, d\^u \`a Yuan~\cite[th.~1.2]{yuan2006}.
\begin{lemm}[Yuan]\label{lemm.yuan}
Soit $\bar L$ et~$\bar M$ des fibr\'es en droites amples sur~$X$
munis de m\'etriques ad\'eliques num\'eriquement
arithm\'etiquement effectifs.
Alors, 
\[ \hvol(\bar L\otimes\bar M^{-1})\geq (\hc_1(\bar L)^{d+1}|X)-(d+1)(\hc_1(\bar L)^d\hc_1(\bar M)|X).\]
\end{lemm}
\begin{proof}
Sous l'hypoth\`ese suppl\'ementaire que
$\bar L$ et $\bar M$ sont munis de m\'etriques arithm\'etiquement
amples,
cette in\'egalit\'e est le th\'eor\`eme principal de~\cite{yuan2006}.
Comme l'hypoth\`ese sur~$\bar L$ et~$\bar M$ entra\^{\i}ne que leurs
m\'etriques ad\'eliques sont limites de m\'etriques arithm\'etiquement amples,
la continuit\'e du volume arithm\'etique implique que
l'in\'egalit\'e est v\'erifi\'ee  dans le cas g\'en\'eral.
\end{proof}

\begin{prop}\label{prop.yuan}
Supposons que $L$ soit ample et que sa m\'etrique soit semi-positive.
Pour tout fibr\'e en droites m\'etris\'e admissible~$\bar M$, on a,
lorsque $t\ra+\infty$,
\begin{align*}
\frac1{t^d} \vol(L^{\otimes t}\otimes M)
   & = (c_1(L)^d|X)+\frac1t d(c_1(L)^{d-1}c_1(M)|X)+\mathrm O(1/t^2) \\
\frac1{t^{d+1}} \hvol(\bar L^{\otimes t}\otimes\bar M)
   &   \geq (\hc_1(\bar L)^{d+1}|X)
      + \frac1t (d+1) (\hc_1(\bar L)^d\hc_1(\bar M)|X)+\mathrm O(1/t^2).
\end{align*}
\end{prop}
\begin{proof}
La premi\`ere \'egalit\'e r\'esulte du th\'eor\`eme de
Hilbert-Samuel, car $L^{\otimes t}\otimes M$ est ample pour $t$
assez grand.

Choisissons des fibr\'es en droites amples
$M_1$ et~$M_2$ sur~$X$  munis de m\'etriques semi-positives
telles que $\bar M\simeq\bar M_1\otimes(\bar M_2)^{-1}$.
Quitte \`a les remplacer par $\bar M_1(c)$ et $\bar M_2(c)$,
o\`u $c$ est un nombre r\'eel assez grand, on suppose qu'ils
sont arithm\'etiquement amples.

Fixons un nombre r\'eel~$c$ tel que $\bar L(c)$
soit arithm\'etiquement num\'eriquement effectif.
Les deux fibr\'es en droites m\'etris\'es $\bar L(c)^t\otimes\bar M_1$
et $\bar M_2$ v\'erifient l'hypoth\`ese du lemme~\ref{lemm.yuan}.
Par suite,
\begin{align*}
 \hvol(\bar L(c)^{\otimes t}\otimes\bar M)
& = \hvol((\bar L(c)^{\otimes t}\otimes\bar M_1)\otimes(\bar M_2)^{-1}) \\
& \geq (\hc_1(\bar L(c)^{\otimes t}\otimes \bar M_1)^{d+1}|X)
    - (d+1) (\hc_1(\bar L(c)^{\otimes t}\otimes \bar M_1)^{d}\hc_1(\bar M_2)|X)\\
& \geq t^{d+1} (\hc_1(\bar L(c))^{d+1}|X)
 + (d+1)t^d (\hc_1(\bar L(c))^d\hc_1(\bar M)|X) + \mathrm O(t^{d-1}) \\
& \geq t^{d+1} (\hc_1(\bar L)^{d+1}|X) +  
(d+1) t^d (\hc_1(\bar L)^d\hc_1(\bar M)|X) \\ 
& \qquad + c(d+1) t^{d+1}\deg_L(X) + c d(d+1) t^d (c_1(L)^{d-1}c_1(M)|X)   \\
& \qquad + \mathrm O(t^{d-1}).
\end{align*}
Par ailleurs,
\[ \hvol(\bar L(c)^{\otimes t}\otimes\bar M)
 = \hvol(\bar L^{\otimes t}\otimes\bar M) + (d+1)t c \vol(L^{\otimes t}\otimes M).\]
Compte-tenu de la premi\`ere \'egalit\'e,
on a donc
\begin{align*}
 \hvol(\bar L^{\otimes t}\otimes\bar M)
-t^{d+1} (\hc_1(\bar L)^{d+1}|X) - (d+1)t^d (\hc_1(\bar L)^d\hc_1(\bar M)|X)
\hskip-20em \\
&\geq c(d+1) t^{d+1}\deg_L(X) + c d(d+1)t^d (c_1(L)^{d-1}c_1(M)|X) \\
& \qquad\qquad  {}   - c (d+1) t \vol(L^{\otimes t}\otimes M) 
                 +\mathrm O(t^{d-1}) \\
&\geq \mathrm O(t^{d-1}),
\end{align*}
ce qui termine la d\'emonstration de la proposition.
\end{proof}

\section{D\'emonstration du th\'eor\`eme d'\'equidistribution logarithmique}
\label{sec.equilog}

Pla\c{c}ons nous sous les hypoth\`eses du th\'eor\`eme~\ref{theo.equilog}.
Fixons en particulier une suite $(x_n)$ de points de $X(\bar F)$,
g\'en\'erique,
telle que $h_{\bar L}(x_n)$ tende vers \mbox{$h_{\bar L}(X)/(d+1)\deg_L(X)$}.

\begin{lemm}\label{lemm.1}
Soit $\bar M$ un fibr\'e en droites muni d'une m\'etrique admissible sur~$X$.
Alors, lorsque $n$ tend vers l'infini,
$ h_{\bar M}(x_n)$ tend vers
\[ \frac{\big(\hc_1(\bar L)^d\hc_1(\bar M)|X \big)}{\deg_L(X)}
    -\frac{d h_{\bar L}(X) \big(c_1(L)^{d-1}c_1(M)|X \big)}{(d+1)\deg_L(X)^2} .\]
\end{lemm}
\begin{proof}
Suivant la voie inaugur\'ee par~\cite{szpiro-u-z97}, la
d\'emonstration repose sur l'existence de sections globales de
normes sup.\ contr\^ol\'ees pour les fibr\'es en droites
m\'etris\'es.

Pour tout nombre entier assez grand~$t$,
le fibr\'e en droites $L^{\otimes t}\otimes M$
est gros. D'apr\`es le lemme~\ref{lemm.inegalite},
\[
\liminf_n \left( t h_{\bar L}(x_n)+h_{\bar M}(x_n)\right) \geq
\frac{\hvol(\bar L^{\otimes t}\otimes \bar M)}{(d+1)\vol(L^{\otimes
t}\otimes M)}.
\]
Comme $h_{\bar L}(x_n)$ converge vers $h_{\bar L}(X)/(d+1)\deg_L(X)$,
le membre de gauche est \'egal \`a
\[ t \frac{h_{\bar L}(X)}{(d+1)\deg_L(X)} + \liminf h_{\bar M}(x_n).\]
D'apr\`es la proposition~\ref{prop.yuan}, le membre de droite
v\'erifie le d\'eveloppement asymptotique suivant lorsque $t\ra\infty$
et $\eps=1/t$:
\begin{align*}
\frac{\hvol(\bar L^{\otimes t}\otimes\bar M)}{\vol(L^{\otimes t}\otimes M)}
& \geq t \frac{(\hc_1(\bar L)^{d+1}|X)+\eps (d+1)(\hc_1(\bar L)^d\hc_1(\bar M)|X)
               +\mathrm O(\eps^2)}
         {(c_1(L)^d|X)+\eps d(c_1(L)^{d-1}c_1(M)|X)+\mathrm O(\eps^2)} \\
& \geq  t \frac{(\hc_1(\bar L)^{d+1}|X)}{(c_1(L)^d|X)}
\left( 1- \eps d \frac{(c_1(L)^{d-1}c_1(M)|X)}{(c_1(L)^d|X)}+\mathrm O(\eps^2)\right)
\\*
& \quad
   +  (d+1) \frac{(\hc_1(\bar L)^d\hc_1(\bar M)|X)}{(c_1(L)^d|X)}
   +\mathrm O(\eps)\\
& = t \frac{(\hc_1(\bar L)^{d+1}|X)}{(c_1(L)^d|X)} \\*
& \quad
{} + (d+1)    \frac{(\hc_1(\bar L)^d\hc_1(\bar M)|X)}{(c_1(L)^d|X)}
 -  d \frac{(\hc_1(\bar L)^{d+1}|X) (c_1(L)^{d-1}c_1(M)|X)}{(c_1(L)^d|X)^2}  \\
& \quad {} +
\mathrm O(\eps).
\end{align*}
Lorsque $t$ tend vers~$+\infty$,  on obtient donc
\[ \liminf_n h_{\bar M}(x_n)
 \geq \frac{(\hc_1(\bar L)^d\hc_1(\bar M)|X)}{\deg_L(X)}
  -   \frac {d h_{\bar L}(X)(c_1(L)^{d-1}c_1(M)|X)}{(d+1)\deg_L(X)^2}
.\]
Lorsqu'on change le fibr\'e en droites m\'etris\'e~$\bar M$ en son inverse,
le membre de gauche
est remplac\'e par $-\limsup_n h_{\bar M}(x_n)$ et le membre
de droite est chang\'e en son oppos\'e,
d'o\`u le lemme.
\end{proof}

Soit $\bar M$ un fibr\'e en droites muni d'une m\'etrique admissible sur~$X$;
soit $s$ une section globale non nulle de~$\bar M$ et soit
$D$ son diviseur. Alors, pour toute place~$v$ de~$F$,
la fonction $\log\norm{s}_v$ est \`a singularit\'es au plus
logarithmiques le long de~$D$.

\begin{lemm}\label{lemm.2}
Pour toute place~$v$ de~$F$,
\[
\liminf_n \int_{X_v} \log\norm{s}^{-1}_v\, \mathrm d\mu_{x_n,v}
\geq \int_{X_v}\log\norm{s}_v^{-1}\,\mathrm d\mu_{\bar L,v}.\]
En outre,
\[
\liminf_n h_{\bar M}(x_n) \geq \sum_v \int_{X_v}\log\norm{s}_v^{-1}\,\mathrm d\mu_{\bar L,v}.
\]
\end{lemm}
\begin{proof}
Fixons une place~$v$ de~$F$.
Pour tout nombre r\'eel~$B$, la fonction $\min(B,\log\norm{s}^{-1}_v)$
est continue sur $X_v$. Le th\'eor\`eme
d'\'equidistribution~\ref{theo.equi} implique donc que,
lorsque $n\ra\infty$,
\[ \int_{X_v} \min(B,\log\norm{s}^{-1}_v)\,\mathrm d\mu_{x_n,v} \ra \int_{X_v}
   \min(B,\log\norm{s}^{-1}_v)\,\mathrm d\mu_{\bar L,v}. \]
Pour tout entier~$n$,
\[ \int_{X_v}\log\norm{s}^{-1}_v\,\mathrm d\mu_{x_n,v}
\geq \int_{X_v}\min( B,\log\norm{s}^{-1}_v)\,\mathrm d\mu_{x_n,v}, \]
si bien que
\begin{multline*}
 \liminf_n \int_{X_v}\log\norm{s}^{-1}_v\,\mathrm d\mu_{x_n,v}
\geq \liminf_n \int_{X_v}\min( B,\log\norm{s}^{-1}_v)\,\mathrm d\mu_{x_n,v} \\
 = \int_{X_v}\min(B,\log\norm{s}^{-1}_v)\,\mathrm d\mu_{\bar L,v} . \end{multline*}
Il reste \`a faire tendre $B$ tend vers~$+\infty$ pour obtenir
la premi\`ere assertion de l'\'enonc\'e.

Montrons maintenant la seconde.
Par d\'efinition,  pour tout entier~$n$ tel que $x_n\not\in D$,
\[ h_{\bar M}(x_n) = \sum_v \int_{X_v}\log\norm{s}_v^{-1}\,\mathrm d\mu_{x_n,v}.\]
La suite~$(x_n)$ \'etant g\'en\'erique, on a donc
\[\liminf_n h_{\bar M}(x_n) \geq \sum_v \liminf_n \int_{X_v}\log\norm{s}_v^{-1}\,\mathrm d\mu_{x_n,v} \geq \sum_v \int_{X_v}\log\norm{s}_v^{-1}\,\mathrm d\mu_{\bar L,v}.\]
Cela termine la d\'emonstration du lemme.
\end{proof}

\begin{proof}[D\'emonstration du th\'eor\`eme~\ref{theo.equilog}]
Soit $M$ le fibr\'e en droites $\mathscr O(D)$; munissons-le
d'une m\'etrique ad\'elique admissible.
Pour toute place~$v$, posons
\[ \eps_v = \liminf_n \int_{X_v}\log\norm s_v^{-1}\,\mathrm d\mu_{x_n,v}
 - \int_{X_v}\log\norm s_v^{-1}\,\mathrm d\mu_{\bar L,v}.\]
En vertu du lemme~\ref{lemm.2}, on a $\eps_v\geq 0$ ; en outre,
\[\liminf_n h_{\bar M}(x_n) \geq \sum_v \int_{X_v}\log\norm{s}_v^{-1}\,\mathrm d\mu_{\bar L,v} + \sum_v \eps_v. \]
Appliquons le lemme~\ref{lemm.1} et le th\'eor\`eme~\ref{theo.mahler},
il vient:
\begin{align*}
\frac1{\deg_L(X)}h_{\bar L}(D)
& = \frac{(\hc_1(\bar L)^d\hc_1(\bar M)|X)}{(c_1(L)^d|X)} - \sum_v\int_{X_v}\log\norm{s}_v^{-1} \mu_{\bar L,v} \\
& =\left( \lim_n h_{\bar M}(x_n)-\sum_v\int_{X_v}\log\norm s_v^{-1}\mu_{\bar L,v}\right) \\
& \qquad {} 
 + \frac {h_{\bar L}(X)}{(d+1)\deg_L(X)} \frac{d \deg_L(D)}{\deg_L(X)} \\
& \geq \sum_v \eps_v
 + \frac {h_{\bar L}(X)}{(d+1)\deg_L(X)} \frac{d \deg_L(D)}{\deg_L(X)} .
\end{align*}
Par cons\'equent,
\[
\frac{h_{\bar L}(D)}{d\deg_L(D)} - \frac{h_{\bar L}(X)}{(d+1)\deg_L(X)}
 \geq  \frac{\deg_L(X)}{d\deg_L(D)} \sum_v \eps_v . \]
Comme $\eps_v\geq 0$ pour tout~$v$, cela d\'emontre
la premi\`ere assertion du th\'eor\`eme.

Supposons que l'on ait l'\'egalit\'e
\[ \frac{h_{\bar L}(D)}{d\deg_L(D)} = \frac{h_{\bar L}(X)}{(d+1)\deg_L(X)}.\]
Alors, pour toute place~$v$, on a $\eps_v=0$, autrement dit
\[ \liminf_n \int_{X_v}\log\norm{s}_v^{-1} \mu_{x_n,v}=\int_{X_v}\log\norm{s}_v^{-1}\mu_{\bar L,v}. \]
Comme toute sous-suite de la suite~$(x_n)$ est encore g\'en\'erique,
on a donc
\[ \lim_n \int_{X_v}\log\norm{s}_v^{-1} \mu_{x_n,v}=\int_{X_v}\log\norm{s}_v^{-1}\mu_{\bar L,v}. \]

L'\'equidistribution pour une fonction \`a singularit\'es
au plus logarithmiques le long de~$D$ en d\'ecoule en vertu
du lemme suivant.
\end{proof}

\begin{lemm}\label{lemm.fatou}
Soit $X$ un espace topologique compact, soit $(\mu_n)$ une suite de
mesures de probabilit\'e sur~$X$ qui converge vers une mesure~$\mu$.
Soit $\phi$ une fonction mesurable positive sur~$X$, $\mu$-int\'egrable
et telle que $\mu_n(\phi)\ra\mu(\phi)$.  
Soit $f$ une fonction mesurable sur~$X$ telle que $\abs{f}\leq \phi$;
si $f$ et $\varphi$ sont continues $\mu$-presque partout, alors
$\mu_n(f)\ra\mu(f)$.
\end{lemm}
\begin{proof}
La fonction $f$ est $\mu$-int\'egrable et on peut supposer qu'elle est
positive. Par d\'efinition m\^eme de la topologie vague sur l'espace
des mesures,  la suite $(\mu_n(g))$ converge vers $\mu(g)$ pour toute
fonction r\'eelle continue $g$ sur $X$ ; plus g\'en\'eralement, ce
r\'esultat vaut pour toute fonction mesurable $g$, born\'ee et continue
$\mu$-presque partout (\cite{bourbaki1965}, Chap.~IV, 5, prop.~22). 
Ceci s'applique en particulier, pour tout nombre r\'eel positif $B$,
aux fonctions $f_B = \min(f,B)$ et $\varphi_B = \min(\varphi,B)$.

Soit $\varepsilon$ un nombre r\'eel strictement positif. La fonction
$\varphi$ \'etant $\mu$-int\'egrable, il existe un nombre r\'eel
positif $B$ tel que $\mu(\varphi - \varphi_B) \leq \varepsilon$. Comme
$\lim_n \mu_n(\varphi) = \mu(\varphi)$ et $\lim_n \mu_n(\varphi_B)
= \mu(\varphi_B)$, la suite $(\mu_n(\varphi - \varphi_B))$ converge
vers $\mu(\varphi - \varphi_B)$ et il existe donc un entier $n_0$ tel
que 
\[ \mu_n(\varphi - \varphi_B) \leq 2\varepsilon\]
pour tout $n \geq n_0$. 
Soit enfin $n_1$ un entier naturel tel que  $\abs{\mu(f_B) - \mu_n(f_B)}
\leq \varepsilon$ pour tout $n \geq n_1$.

En utilisant l'encadrement 
\[ 0 \leq f - f_B = \max(f,B) - B \leq \max(\varphi,B) - B 
= \varphi - \varphi_B,\]
on d\'eduit de ce qui
pr\'ec\`ede la majoration 
\[ \abs{\mu(f) - \mu_n(f)} \leq \mu(\varphi - \varphi_B) + \abs{\mu(f_B) - \mu_n(f_B)} + \mu_n(\varphi - \varphi_B) \leq
4\varepsilon\]
pour tout entier $n \geq \max(n_0,n_1)$ et la suite
$(\mu_n(f))$ converge bien vers $\mu(f)$.
\end{proof}

\section{\'Equidistribution logarithmique de sous-vari\'et\'es}

Nous d\'emontrons ici le th\'eor\`eme~\ref{theo.equilog.ssvar}.
Rappelons que l'on pose, pour toute sous-vari\'et\'e~$Y$ de~$X$,
\[ h'_{\bar L}(Y)= \frac{(\hc_1(\bar L)^{1+\dim Y}|Y)}{(1+\dim(Y)) (c_1(L)^{\dim Y}|Y)}. \]
Si $p$ est un entier, notons $\mu_p(\bar L)$ la borne inf\'erieure
des $h'_{\bar L}(Y)$, lorsque $Y$ parcourt l'ensemble
des sous-vari\'et\'es de dimension~$p$ de~$X$.
Si $L$ est ample, alors $\mu_p(\bar L)$ est un nombre r\'eel
(\cite{bost-g-s94}, remarque~(iii) p.~954).
Notons aussi $\lambda_p(\bar L)$ la quantit\'e suivante:
\[ \lambda_p(\bar L) = \sup_{\substack{\varnothing \neq U\subset X \\ \text{$U$ ouvert}}}
 \inf_{\substack{Y\cap U\neq\varnothing \\ \dim Y=p}} h'_{\bar L}(Y)\quad;\]
en particulier, $e(\bar L)=\lambda_0(\bar L)$.

\begin{lemm}\label{lemm.1.ssvar}
Supposons que $\bar L$ soit un fibr\'e en droites ample muni d'une
m\'etrique ad\'elique semi-positive.
Soit $p$ un entier naturel et
soit $(Y_n)$ une suite g\'en\'erique de sous-vari\'et\'es de~$X$
de dimension~$p$
telle que $h'_{\bar L}(Y_n)$ tende vers $h'_{\bar L}(X)$.
Soit aussi $\bar M$ un fibr\'e en droites admissible~$\bar M$ sur~$X$.

Supposons que l'on ait $\mu_{p-1}(\bar L)\geq h'_{\bar L}(X)$.
Alors, lorsque $n$ tend vers l'infini, 
\[ \frac{ (\hc_1(\bar L)^p \hc_1(\bar M)|Y_n)}{(c_1(L)^p|Y_n)}
 - p  \frac{(c_1(L)^{p-1}c_1(M)|Y_n)}{(c_1(L)^p|Y_n)} h'_{\bar L}(X) \]
converge vers
\[  \frac{(\hc_1(\bar L)^d \hc_1(\bar M)|X)}{(c_1(L)^d|X)}
 - d h'_{\bar L}(X) \frac{(c_1(L)^{d-1}c_1(M)|X)}{(c_1(L)^d|X)}. \]
\end{lemm}
\begin{proof}
Cela se d\'emontre de mani\`ere similaire au lemme~\ref{lemm.1}.
Toutefois, pour simplifier les calculs, faisons l'hypoth\`ese suppl\'ementaire
que $h'_{\bar L}(X)=0$ quitte \`a multiplier la m\'etrique
de~$\bar L$ par une constante aux places archim\'ediennes.
On v\'erifie facilement que l'\'enonc\'e obtenu sous cette hypoth\`ese
est \'equivalent au r\'esultat que nous voulons d\'emontrer.
 
Lorsque $t$ tend vers l'infini, $L^{\otimes t}\otimes M$ est gros
et la proposition~\ref{prop.yuan} entra\^{\i}ne la minoration
\[ \frac{\hvol(\bar L^{\otimes t}\otimes\bar M)}{(d+1)\vol(L^{\otimes t}\otimes M)}
 \geq \frac{(\hc_1(\bar L)^d \hc_1(\bar M)|X)}{(c_1(L)^d|X)} + \mathrm O(1/t).
\]
Posons $\bar L_t=\bar L^{\otimes t}\otimes \bar M$;
soit $\eps$ un nombre r\'eel compris entre~$0$ et~$1$,
soit $k$ un entier positif tel que
\[ (1-\eps) \vol(L_t) \frac{k^d}{d!} \leq \dim H^0(X,L_t^{\otimes k})\leq (1+\eps) \vol(L_t) \frac{k^d}{d!} \]
et
\[ \chi_{\sup}(\bar L_t^{\otimes k}) \geq (\hvol(\bar L_t)-\eps) \frac{k^{d+1}}{(d+1)!}. \]
Posons $D=h^0(X,L_t^{\otimes k})$ et soit
$c$ le nombre r\'eel tel que $c^D \exp(\chi_{\sup}(\bar L_t^{\otimes k})) = 2^D$.
Soit alors $s$ une section non nulle de $\bar L_t^{\otimes k}$
qui, 
comme dans la d\'emonstration du lemme~\ref{lemm.inegalite},
v\'erifie les in\'egalit\'es $\norm s_v\leq 1$  pour $v$ ultram\'etrique
et $\norm s_v\leq c$ pour $v$ archim\'edienne.
Comme la suite~$(Y_n)$ est g\'en\'erique, $Y_n\not\subset \div(s)$
pour $n$ assez grand, c'est-\`a-dire $s_{|Y_n}\neq 0$.
Alors,
\begin{align*}
(\hc_1(\bar L)^p\hc_1(\bar M)|Y_n) \hskip -4em \\
&=
 (\hc_1(\bar L)^p \hc_1( \bar L^{\otimes t}\otimes \bar M)|Y_n)
 - (\hc_1(\bar L)^{p+1}|Y_n) \\
&= \frac1k (\hc_1(\bar L)^p | \div(s_{|Y_n}))
 + \sum_v \int_{X_v} \log\norm{s}_v^{-1/k} c_1(\bar L)^p\delta_{Y_n}
 - (\hc_1(\bar L)^{p+1}|Y_n) \\
& \geq \frac 1k (\hc_1(\bar L)^p | \div(s_{|Y_n}))
- (c_1(L)^p|Y_n)  \frac1k \log c
- (p+1) (c_1(L)^p|Y_n) h'_{\bar L}(Y_n).
\end{align*}
Puisque $\div(s_{|Y_n})$ est un cycle effectif de dimension~$p-1$,
on a 
\[ (\hc_1(\bar L)^p | \div(s_{|Y_n})) \geq p (c_1(L)^{p-1}|\div(s_{Y_n})) 
h'_{\bar L}(X)
 \geq 0. \]
Par suite, lorsque $n$ tend vers l'infini,
\begin{align*}
\liminf \frac{(\hc_1(\bar L)^p\hc_1(\bar M)|Y_n)}{(c_1(L)^p|Y_n)}
&  \geq  - \frac1k \log c - (p+1) h'_{\bar L}(X) 
 \geq -\frac 1k \log c \\
& \geq  \frac1{Dk} \chi_{\sup}(\bar L_t^{\otimes k})
 - \frac 1k \log 2 \\
& \geq \min \big(\frac{\hvol(\bar L_t)-\eps}{(d+1)\vol(L_t)(1-\eps)},
  \frac{\hvol(\bar L_t)-\eps}{(d+1)\vol(L_t)(1+\eps)} \big). \end{align*}
Par suite,
\[ \liminf \frac{(\hc_1(\bar L)^p\hc_1(\bar M)|Y_n)}{(c_1(L)^p|Y_n)}
\geq \frac{\hvol(\bar L _t)}{(d+1)\vol(L_t)}
\geq  
\frac{(\hc_1(\bar L)^d \hc_1(\bar M)|X)}{(c_1(L)^d|X)} + \mathrm O(1/t).
\]
Faisons maitenant tendre $t$ vers l'infini; on en d\'eduit
\[  \liminf \frac{(\hc_1(\bar L)^p\hc_1(\bar M)|Y_n)}{(c_1(L)^p|Y_n)}
\geq \frac{(\hc_1(\bar L)^d \hc_1(\bar M)|X)}{(c_1(L)^d|X)} .\]

L'in\'egalit\'e dans l'autre sens en d\'ecoule alors
en rempla\c{c}ant $\bar M$ par son inverse.
\end{proof}

\begin{rema}
Supposons que $\bar M$ soit le fibr\'e en droite trivial 
muni de la m\'etrique triviale \`a toute place, sauf en une place~$v$ de~$F$;
posons $\phi_v=-\log\norm{1}_v$.
Sous les hypoth\`eses du lemme~\ref{lemm.1.ssvar}, on a
\[ \lim_{n\ra+\infty}
\int_{X_v} \phi_v  \frac{c_1(L)^p_v\delta_{Y_n}}{(c_1(L)^p|Y_n)}
= \int_{X_v} \phi_v \mathrm d\mu_{\bar L,v}. \]
Comme les fonctions admissibles sont denses dans les fonctions
continues, on retrouve ainsi le th\'eor\`eme d'\'equidistribution des sous-vari\'et\'es
qu'ont d\'emontr\'e Autissier~\cite{autissier2006b} et Yuan~\cite{yuan2006}.
\end{rema}

Pla\c{c}ons-nous sous les hypoth\`eses du th\'eor\`eme~\ref{theo.equilog.ssvar}.
Pour all\'eger les notations, on notera  pour tout entier~$n$ et toute
place~$v$ de~$F$.
\[ \mathrm d\mu_{Y_n,v} = \frac{1}{(c_1(L)^p|Y_n)} c_1(\bar L)^p_v \ ; \]
c'est une mesure sur~$X_v$ support\'ee par $Y_n$. 
D'apr\`es la remarque pr\'ec\'edente, pour toute place~$v$
de~$F$, $\mathrm d\mu_{Y_n,v}\ra \mathrm d\mu_{\bar L,v}$
lorsque $n$ tend vers~$+\infty$. Il s'agit de d\'emontrer
que la convergence a encore lieu lorsque ces mesures sont
multipli\'ees par certaines fonctions \`a singularit\'es au plus logarithmiques.

Supposons de plus que $h'_{\bar L}(X)=0$ ; on laisse au lecteur
le soin de v\'erifier que cela ne restreint pas la g\'en\'eralit\'e.

Soit $\bar M$ un fibr\'e en droites muni d'une m\'etrique admissible
sur~$X$, soit $s$ une section globale non nulle de~$\bar M$
et soit $D$ son diviseur. 
Pour toute place~$v$ de~$F$, la fonction $\log\norm s_v$ est \`a singularit\'es
au plus logarithmiques le long de~$D$.

\begin{lemm}
Pour toute place~$v$ de~$F$,
\[ \liminf_n \frac1{(c_1(L)^p|Y_n)}
 \int_{X_v} \log\norm s_v^{-1} c_1(\bar L)^p_v \delta_{Y_n}
 \geq \int_{X_v} \log \norm{s}_v^{-1} \,\mathrm d\mu_{\bar L,v} .\]
De plus,
\[ \liminf \frac{(c_1(\bar L)^pc_1(\bar M)|Y_n)}{(c_1(L)^p|Y_n)}
 \geq \sum_v \int_{X_v}  \log \norm{s}_v^{-1} \,\mathrm d\mu_{\bar L,v} .\]
\end{lemm}
\begin{proof}
Il s'agit d'une simple adaptation du lemme~\ref{lemm.2}.
La premi\`ere in\'egalit\'e se d\'emontre en appliquant la remarque pr\'ec\'edente \`a 
la fonction $\min(B,\log\norm{s}_v)$ sur~$X_v$
et en faisant tendre $B$ vers~$+\infty$.
Pour la seconde, on rappelle (th\'eor\`eme~\ref{theo.mahler}) que
\[ (\hc_1(\bar L)^p\hc_1(\bar M)|Y_n) = \sum_v \int_{X_v} \log\norm{s}_v^{-1}\, c_1(\bar L)^p_v\delta_{Y_n} + (\hc_1(\bar L)^p| \div(s_{|Y_n})). \]
Par cons\'equent,
\begin{align*}
  \liminf \frac{(c_1(\bar L)^pc_1(\bar M)|Y_n)}{(c_1(L)^p|Y_n)}
 & \geq \sum_v \liminf_n \int_{X_v} 
 \log\norm s_v^{-1} c_1(\bar L)^p_v \delta_{Y_n} \\
&  \geq \sum_v \int_{X_v} \log \norm{s}_v^{-1} \,\mathrm d\mu_{\bar L,v}.
\qquad\qedhere\end{align*}
\end{proof}

Pour toute place~$v$ de~$F$, posons
\[ \eps_v =\liminf_n  \int_{X_v} \log\norm{s}_v^{-1} \,\mathrm d\mu_{Y_n,v}
 - \int_{X_v} \log\norm{s}_v^{-1} \mathrm d\mu_{\bar L,v}. \]
On vient de d\'emontrer que $\eps_v\geq 0$ ; en outre,
\[ \liminf \frac{(c_1(\bar L)^pc_1(\bar M)|Y_n)}{(c_1(L)^p|Y_n)}
\geq \sum_v \int_{X_v} \log\norm{s}_v^{-1} \,\mathrm d\mu_{\bar L,v}
 + \sum_v \eps_v. \]
Appliquons maintenant le lemme~\ref{lemm.1.ssvar} ; on a ainsi
\begin{align*}
\frac1{(c_1(L)^d|X)} h_{\bar L}(D)
& =  \frac{(\hc_1(\bar L)^d\hc_1(\bar M)|X)}{(c_1(L)^d|X)}
- \sum _v \int_{X_v} \log\norm s_v^{-1} \,\mathrm d\mu_{\bar L,v} \\
& =  \lim_n \frac{(\hc_1(\bar L)^p\hc_1(\bar M)|X)}{(c_1(L)^d|X)}
- \sum _v \int_{X_v} \log\norm s_v^{-1} \,\mathrm d\mu_{\bar L,v} 
\\
&\geq \sum_v \eps_v. 
\end{align*}

Puisqu'on a suppos\'e  $h_{\bar L}(D)=0$, il vient $\eps_v=0$ pour tout~$v$.
Autrement dit,
\[ \liminf_n \int_{X_v} \log\norm s_v^{-1} \,\mathrm d\mu_{Y_n,v}
 = \int_{X_v} \log\norm s_v^{-1} \,\mathrm d\mu_{\bar L,v}. \]
Le r\'esultat analogue pour la limite, plut\^ot qu'une limite inf\'erieure,
s'en d\'eduit en consid\'erant une sous-suite de la suite~$(Y_n)$.
Le th\'eor\`eme~\ref{theo.equilog.ssvar}
pour une fonction \`a singularit\'es au plus logarithmiques
le long de~$D$  g\'en\'erale en d\'ecoule via le lemme~\ref{lemm.fatou}.

Observons enfin que,
comme dans le cas des points, 
les hypoth\`eses du th\'eor\`eme~\ref{theo.equilog.ssvar}  
impliquent que $h'_{\bar L}(D)\geq h'_{\bar L}(X)$ 
pour tout diviseur effectif~$D$. 

\section{Non-\'equidistribution logarithmique de certaines suites de points}
\label{sec.autissier}

L'\'enonc\'e suivant nous a \'et\'e sugg\'er\'e par P.~Autissier.
\begin{prop}
Soit $X$ un sch\'ema projectif et int\`egre de dimension $d$. Soit
$\bar L$ un fibr\'e en droites sur $X$ muni d'une m\'etrique
ad\'elique semi-positive.

Soit $\bar M$ un fibr\'e en droites sur~$X$ muni d'une m\'etrique
admissible; soit $s$ une section globale non nulle de~$M$ et soit
$D$ son diviseur.

Soit $(x_{n})$ une suite g\'en\'erique de points ferm\'es de $X$
telle que $h_{\bar L}(x_{n})$ converge vers $h_{\bar
L}(X)/(d+1)\deg_L(X)$.

Fixons un ensemble fini $S$ de places de $F$ contenant toutes les
places archim\'ediennes et supposons en outre qu'il existe un
$\mathscr O_{F}$-sch\'ema projectif et plat $\mathscr X$,
de fibre g\'en\'erique $X$,
ainsi qu'un fibr\'e en droites $\mathscr L$
et un diviseur de Cartier horizontal et effectif~$\mathscr D$
sur $\mathscr X$ tels que :
\begin{itemize}
\item $\mathscr L_{|X} \simeq L$
et,
en toute place~$v\not\in S$,
la m\'etrique de $\bar L$ est lisse, d\'efinie
par ce mod\`ele~$(\mathscr X,\mathscr L)$ ;
\item $\mathscr D_{|X} = D$ et la m\'etrique de~$\bar M$
est d\'efinie par $\mathscr D$ en toute place~$v\not\in S$;
\item les points~$x_{n}$ se prolongent en des points $S$-entiers
de $\mathscr X \setminus \mathscr D$.
\end{itemize}

Dans ces conditions, la suite $(\sum_{v \in S} \int_{X_v} \log\norm{s}^{-1}_v
\mathrm d\mu_{x_n, v} ) $
converge vers
\[
\sum_{v \in S} \int_{X_v} \log\norm{s}^{-1}_v
\mathrm d\mu_{\bar L, v} + \frac{d \deg_L (D)}{\deg_L
(X)}\left(\frac{h_{\bar L}(D)}{d \deg_L (D)} - \frac{h_{\bar
L}(X)}{(d+1) \deg_L (X)} \right).
\]
\end{prop}

Lorsque l'hypoth\`ese du th\'eor\`eme
\ref{theo.equilog} sur la hauteur du diviseur $D$ n'est pas
v\'erifi\'ee,
cela montre que pour au moins une place $v\in S$,
il n'y a pas \'equidistribution pour la fonction
$\log\norm{s}_v^{-1}$.

\begin{proof}
Soit $v$ une place de $F$ n'appartenant pas \`a $S$.
Le diviseur~$\mathscr D$ est effectif et d\'efinit
la fonction de Green~$\log\norm s_v^{-1}$
pour~$D$; en particulier, $\log\norm s_v^{-1}\geq 0$
et pour un point $x\in X(\C_v)$,
le cas d'\'egalit\'e $\log\norm {s(x)}_v^{-1}=0$
\'equivaut \`a ce que le point~$x$ se prolonge
en un point de~$\mathscr X\setminus\mathscr D$ au-dessus
de~$v$.
Par hypoth\`ese, $x_{n}$ s'\'etend en un point
$S$-entier de $\mathscr X \setminus \mathscr D$ ;
on a donc
$\int_{X_v} \log\norm{s}_{v}^{-1} \mathrm d\mu_{x_n,v} = 0$
pour toute place~$v$ de~$F$ telle que $v\not\in S$.

Comme la m\'etrique de $\bar L$ \`a la place~$v$ est lisse et
provient du prolongement~$\mathscr L$ de $L$ sur~$\mathscr X$, la
mesure~$\mathrm d\mu_{\bar L,v}$ est une combinaison lin\'eaire de
masses de Dirac aux points de~$X_v$ qui se r\'eduisent sur les
points g\'en\'eriques de la fibre sp\'eciale~$\mathscr X_v$
de~$\mathscr X$ en~$v$. En un tel point~$\xi$, correspondant \`a une
composante irr\'eductible~$C$, $\log\norm{s(\xi)}^{-1}$ est la
multiplicit\'e de  la composante~$C$ dans~$\mathscr D$
(multipli\'ee par $\log\abs\pi^{-1}$, o\`u $\pi$ est une uniformisante
de~$F$ en~$v$), donc est nul. 
En particulier, $\int_{X_v} \log\norm{s}^{-1}_v \mathrm
d\mu_{\bar L, v} = 0$. Pour tout~$n$, on a donc l'\'egalit\'e
\[ h_{\bar M}(x_n)
= \sum_{v\in S} \int_{X_v} \log\norm s_v^{-1}\,\mathrm d\mu_{x_n,v}. \]
De plus, d'apr\`es le lemme~\ref{lemm.1},
la suite $(h_{\bar M}(x_n))$ converge vers
\[\frac{\big( \hc_1(\bar L)^d\hc_1(\bar M)|X \big)}{\deg_L (X)} -
\frac{d h_{\bar L}(X)\big(c_1(L)^{d-1}c_1(M)|X
\big)}{(d+1)\deg_L(X)^2}. \] 
Compte tenu du th\'eor\`eme~\ref{theo.mahler} et
de la nullit\'e des int\'egrales $\int_{X_v}\log\norm
s_v^{-1}\mathrm d\mu_{\bar L,v}$ si $v\not\in S$, 
on a 
\[\frac{\big(\hc_1(\bar L)^d\hc_1(\bar
M)|X\big)}{\deg_L(X)} = \frac{h_{\bar L}(D)}{\deg_L(X)} + \sum_{v
\in S} \int_{X_v} \log\norm{s}_v^{-1} \mathrm d\mu_{\bar L, v}.
\]
Par suite, $h_{\bar M}(x_n)$ tend
vers \[\sum_{v \in S} \int_{X_v}
\log \norm{s}^{-1}_v \mathrm d\mu_{\bar L, v} + d\frac{\deg_L
(D)}{\deg_L (X)} \left(\frac{h_{\bar L}(D)}{d \deg_L (D)} -
\frac{h_{\bar L}(X)}{(d+1)\deg_L (X)}\right),\] 
ce qui d\'emontre la proposition.
\end{proof}

On peut ici retrouver l'exemple de l'article~\cite{autissier2006}.
Prenant pour $X$ la droite projective ${\P}^1_{\Q} = \Proj (\Q [T_0,
T_1])$, Autissier consid\`ere en effet la suite des points ferm\'es
$x_{n}$ d\'efinis par les polyn\^omes irr\'eductibles $P_n = (T^n -
1)(T - 2)+ 3 \in \Z [T]$ ($T = T_{0}^{-1}T_1$ est la coordonn\'ee
affine de~$\mathbf A^1_\Q$ ). Soit $\bar L$ le fibr\'e tautologique
$\mathscr O(1)$ muni de la m\'etrique de Weil; la suite $(x_n)$ est
g\'en\'erique et la suite des hauteurs $(h_{\bar L}(x_n))$ converge
vers $0=h_{\bar L}(X)$. Cependant, notant $f_{1} $ la fonction $
\log \max(1,\abs{T-2}_{\infty}) = \log\norm{T-2}_{\infty}^{-1}$ sur
$X_{\infty}$, Autissier montre que $\int_{X_{\infty}} f_{1} \mathrm
d \mu_{x_n ,\infty}$ converge vers $\log 2$ et non vers
$\int_{X_{\infty}} f_1 c_1(\bar L)_{\infty} = 0$.

Posons $\mathscr X = \Proj (\Z [T_0, T_1])$
et soit $\mathscr D$ le lieu des z\'eros de~$T-2$.
Prenons pour ensemble~$S$ de places de~$\Q$  la partie $\{\infty, 3\}$.
On voit sur l'expression du polyn\^ome~$P_n$ que  $x_n$
se prolonge en un point $S$-entier de ${\mathscr X} \setminus \mathscr D$.
D'apr\`es la proposition pr\'ec\'edente, la suite
\[ \int_{X_\infty} \log\norm{T-2}_{\infty}^{-1} \mathrm d\mu_{x_n
,\infty} + \int_{X_3} \log\norm{T-2}_{3}^{-1} \mathrm d\mu_{x_n , 3}\]
converge donc vers $h_{\bar L}(2) = \log 2$.
En consid\'erant le polygone de Newton $3$-adique de~$P_n(T-2)$, 
on d\'emontre que pour tout entier~$n$,
\[\int_{X_3} \log\norm{T-2}_{3}^{-1}
\mathrm d\mu_{x_n, 3} = \frac{\log 3}{n+1}. \]
Par cons\'equent, $\int_{X_\infty} \log\norm{T-2}^{-1}_{\infty} \mathrm d\mu_{x_n,
\infty}$ converge vers $\log 2$.

\bibliographystyle{smfplain}
\bibliography{aclab,acl,equilog}
\end{document}